\documentclass[review,3p]{elsarticle}

\usepackage{2_share/style_latex}

\begin{document}

\begin{frontmatter}

\title{A practical algorithm to minimize the overall error in FEM computations}

\author[1]{Jie Liu\corref{cor1}}
\ead{j.liu-5@tudelft.nl}
\author[1]{Henk M. Schuttelaars}
\ead{h.m.schuttelaars@tudelft.nl}
\author[1]{Matthias M\"oller}
\ead{m.moller@tudelft.nl}

\address[1]{Delft Institute of Applied Mathematics\\ Delft University of Technology\\ Mekelweg 4, 2628 CD Delft, The Netherlands}
\cortext[cor1]{Corresponding author}

\begin{abstract}
Using the standard finite element method (FEM) to solve general partial differential equations, the round-off error is found to be proportional to $N^{\beta_{\rm R}}$, with $N$ the number of degrees of freedom (DoFs) and $\beta_{\rm R}$ a coefficient.
A method which uses a few cheap numerical experiments is proposed to determine the coefficient of proportionality and $\beta_{\rm R}$ in various space dimensions and FEM packages.
Using the coefficients obtained above, the strategy put forward in \cite{liu386balancing} for predicting the highest achievable accuracy $E_{\rm min}$ and the associated optimal number of DoFs $N_{\rm opt}$ for specific problems is extended to general problems.
This strategy allows predicting $E_{\rm min}$ accurately for general problems, with the CPU time for obtaining the solution with the highest accuracy $E_{\rm min}$ typically reduced by 60\%--90\%.
\end{abstract}



\end{frontmatter}

\section{Introduction}

Many problems in scientific computing consist of solving boundary value problems.
In this paper, the use of the standard finite element method (FEM) is considered.
To get accurate solutions, various approaches, such as $h$-refinement, $p$-refinement, or $hp$-refinement, are used.
The $h$-refinement is investigated in detail in this paper.

The $h$-refinement typically focuses on the reduction of the truncation error $E_{\rm T}$ by decreasing the grid size, denoted by $h$, of the discretized problem.
However, the round-off error $E_{\rm R}$ increases with the decreasing grid size and will exceed the truncation error when the grid size is getting too small~\cite{Babuska2018Roundoff,butcher2016numerical}.
The dependency of the truncation error $E_{\rm T}$ and the round-off error $E_{\rm R}$ on the grid size results in the total error $E_h = (E_{\rm T} + E_{\rm R})$ first decreasing and then increasing with the decreasing grid size.

For the same grid size, the numerical accuracy might be different for different problems~\cite{liu386balancing}.
To reach the accuracy required, obtained with an element size denoted by $h_{\rm tol}$, the common method is to refine the mesh sequentially from a low refinement level until the required accuracy is satisfied.
This process may take a large number of $h$-refinements.
Since the results for the grid size larger than $h_{\rm tol}$ are thrown away after the required accuracy is satisfied, and the required accuracy may not be reached due to various reasons such as differentiation~\cite{wei2017three}, we call this method the brute-force method, denoted by BF.

To know if the required accuracy can be reached, and if it can be reached in an efficient way, estimates of the highest achievable numerical accuracy, denoted by $E_{\rm min}$, and the associated optimal number of DoFs, denoted by $N_{\rm opt}$, must be available.
In \cite{liu386balancing}, they were predicted using the relations between the round-off error and the number of DoFs and the truncation error and the number of DoFs.
This approach was applied to one-dimensional problems.
It was assumed the behaviour of the round-off error, represented by two coefficients, was known.
This strategy allowed us to predict the accuracy $E_{\rm min}$ using a few computations on coarse grids, of which the CPU time taken is negligible.
Furthermore, for obtaining the solution with the highest accuracy by computing the result using the $N_{\rm opt}$ predicted, the CPU time reduction is around 70\%.

Since only a few specific 1D cases were considered to obtain the coefficients of the round-off error, these coefficients may not apply to other 1D cases.
Their applicability may become even worse when using them directly for 2D problems since the method of counting the number of DoFs and the magnitude of the number of DoFs for 2D problems is different from that for 1D problems.
In view of the above, the aim of this paper is to investigate the coefficients of the round-off error when solving generic 1D and 2D partial differential equations, and extend the strategy to obtain $N_{\rm opt}$, as proposed in \cite{liu386balancing}, to generic 1D and 2D problems. 

The paper is organized as follows. 
The model problem, finite element method, and highest achievable accuracy are discussed in Section \ref{section_model_problem_FE_method_highest_achievable_accuracy}.
A novel method to obtain the coefficients of the round-off error is illustrated in Section \ref{section_novel_approach_for_obtaining_round_off_error}.
The algorithm for determining the above coefficients and predicting the accuracy $E_{\rm min}$ is put forward in Section \ref{section_algorithm}, followed by a validation of our method in Section \ref{section_validation_PRED_plus}.
Conclusions are drawn in Section \ref{paragraph_on_conclusion}.

\captionsetup[figure]{labelfont={bf},name={Fig.},labelsep=period}                                                   
\captionsetup[table]{labelfont={bf},name={Table},labelsep=space,justification=justified,singlelinecheck=false}      

\section{Model problem, finite element method, and highest achievable accuracy}	\label{section_model_problem_FE_method_highest_achievable_accuracy}

\subsection{Model problem}				\label{section_model_problem}

We consider the following second-order partial differential equation:
\begin{equation}
 -\nabla \cdot \left(D(\mathbf{x}) \nabla u \right) + r(\mathbf{x})u(\mathbf{x}) = f(\mathbf{x}),\qquad \mathbf{x} \in \Omega = [0,\,1] \times [0,\,1],	\label{2d_second_order_differential_equation}
\end{equation}
where $u$ denotes the unknown dependent variable, $f(\mathbf{x})$ the prescribed right-hand side, $D(\mathbf{x})$ the diffusion matrix, and $r(\mathbf{x})$ the coefficient function of the reactive term. 
$f(\mathbf{x})$, $D(\mathbf{x})$, and $r(\mathbf{x})$ are continuous and elements of the function space $L_2(\Omega)$. 

The matrix $D(\mathbf{x})$ is considered to be symmetric and positive definite.
By choosing $D(\mathbf{x}) = I$, the identity matrix, and $r(\mathbf{x})=0$, Eq.~(\ref{2d_second_order_differential_equation}) reduces to the Poisson equation;
the diffusion equation is found if $r(\mathbf{x})=0$ and $D(\mathbf{x})$ \enquote{arbitrary} with the above constraints, and the Helmholtz equation is found for $r(\mathbf{x}) \neq $ 0 and $D(\mathbf{x}) = I$.
If not stated otherwise, at the left and right boundaries, denoted by $\Gamma_D$, Dirichlet boundary conditions are imposed: $u(\mathbf{x})=g(\mathbf{x})$. At the upper and bottom boundaries, denoted by $\Gamma_N$, Neumann boundary conditions are prescribed: $D(\mathbf{x})\nabla u \cdot \mathbf{n} = h(\mathbf{x})$, where $\mathbf{n}$ denotes the outward pointing unit normal vector.

\subsection{Finite element method} 	                \label{finite_element_method}

\subsubsection{Weak form}
To construct the finite element (FE) approximation of this equation, we first derive the weak form and then discretize it using proper FE spaces.
For convenience, we introduce three inner products~\cite{lipschutz2009linear}:
\begin{subequations}
  \begin{align}
   \langle \mathbf{f}_1, \,\mathbf{f}_2 \rangle &= \int _{\Omega} \mathbf{f}_1(\mathbf{x}) \cdot \mathbf{f}_2(\mathbf{x}) \, d A,	\\
   \langle f_1, \,f_2 \rangle &= \int _{\Omega} f_1(\mathbf{x}) f_2(\mathbf{x}) \, d A,	\text{ and}\\
   \langle f_1, \,f_2 \rangle _{\Gamma} &= \int _{\Gamma} f_1(\mathbf{x}) f_2(\mathbf{x}) ds,
  \end{align}
\end{subequations}
where $\mathbf{f}_1(\mathbf{x})$ and $\mathbf{f}_2(\mathbf{x})$ denote continuous two-dimensional vector-valued functions, and $f_1(\mathbf{x})$ and $f_2(\mathbf{x})$ denote continuous scalar functions. $\Gamma$ denotes the boundary of $\Omega$.
Moreover, the following function spaces are defined~\cite{chen2005finite}:
\begin{subequations}
\begin{align}
  V_g &= \{t \; | \; t \in H^1 (\Omega), \; t = g \text{ on } \Gamma _D \},  \\
  V_0 &= \{t \; | \; t \in H^1 (\Omega), \; t = 0 \text{ on } \Gamma _D\}.
\end{align}
\label{function_spaces_for_the_weak_form}%
\end{subequations}

\paragraph{Weak form derivation}

Multiplying Eq. (\ref{2d_second_order_differential_equation}) by a test function $\eta \in V_0$~\cite{chen2005finite,gockenbach2006understanding} and integrating it over $\Omega$ yield
\begin{equation}
\langle \eta, \, - \nabla \cdot \left( D \nabla u \right) + ru \rangle = \langle \eta, \, f \rangle. \label{2D_general_inte}
\end{equation}
By applying Gauss's theorem and substituting $\eta = 0$ at $\Gamma _D$, we obtain
\begin{equation}
 \langle {\nabla \eta}, \, D \nabla u \rangle + \langle \eta, \, ru \rangle = \langle \eta, \, f \rangle + \langle \eta, \, D \nabla u \cdot \mathbf{n} \rangle_{ {\Gamma_N}}.	\label{2D_general_gauss}
\end{equation}
Substituting the natural boundary condition, i.e. $D \nabla u \cdot \mathbf{n} = h$ on $\Gamma_N$, we obtain
\begin{equation}
 \langle {\nabla \eta}, \, D \nabla u \rangle + \langle \eta, \, ru \rangle = \langle \eta, \, f \rangle + \langle \eta, \, h \rangle_{ {\Gamma_N}}.		\label{2D_general_neum}
\end{equation}
As a result, the weak form of Eq. (\ref{2d_second_order_differential_equation}) reads
\begin{equation}
\centering
\boxed{
\begin{aligned}
&\text{Find $u \in V_g$ such that:} \\
& \langle {\nabla \eta}, \, D \nabla u \rangle + \langle \eta, \, ru \rangle = \langle \eta, \, f \rangle + \langle \eta, \, h \rangle_{ {\Gamma_N}} \qquad \forall \eta \in V_0.\\
\end{aligned}		\label{sm_weak_form_finite_element_not_specified} 
}
\end{equation}

\noindent The terms on the right-hand side of Eq.~(\ref{sm_weak_form_finite_element_not_specified}) consist of the weakly imposed body force $f$ and the Neumann boundary conditions.
The latter vanishes if no Neumann boundary conditions are prescribed.

\paragraph{Weak form discretization}
This section contains two parts: defining the FE space and constructing the system of equations on the above FE space.
The FE space for $V_g$ reads~\cite{chen2005finite}
\begin{subequations}
\begin{align}
  V_{h,g} = \{ t \; | \; t \text{ is continuous on } \Omega,~ t \in T_p(K), K \in K_h \text{, and } t = g \text{ on } \Gamma _D \},
\end{align}
and that for $V_0$ reads
\begin{align}
  V_{h,0} = \{ t \; | \; t \text{ is continuous on } \Omega,~ t \in T_p(K), K \in K_h \text{, and } t = 0 \text{ on } \Gamma _D \},
\end{align}
\label{finite_element_space_definition}%
\end{subequations}
where $T_p$ denotes functions built by the Lagrangian polynomials of degree $p$, $K$ each individual mesh element, and $K_h$ the computational mesh for $\Omega$.

Using the above FE space, the numerical solution is approximated by
\begin{equation}
 u_h = \sum _ {i=1} ^{m} u _{i} \varphi _i, \label{sm_u_approx}%
\end{equation}
where $\varphi _i$ denotes the basis functions, $m$ is the number of basis functions, and $u_i$ are solution values at the support points of the basis functions, i.e. DoFs.

Consequently, with Eq.~(\ref{finite_element_space_definition}) and Eq.~(\ref{sm_u_approx}), the weak form Eq.~(\ref{sm_weak_form_finite_element_not_specified}) can be discretized as
\begin{equation}
\centering
\boxed{
\begin{aligned}
&\text{Find $u_h \in V_{h,g}$ such that:} \\
& \langle {\nabla \eta}, \, D \nabla u_h \rangle + \langle \eta, \, r u_h \rangle = \langle \eta, \, f \rangle + \langle \eta, \, h \rangle_{ {\Gamma_N}} \qquad \forall \eta \in V_{h,0}.
\end{aligned}
}                           \label{sm_weak_form_discretized}
\end{equation}

\noindent Choosing the test functions $\eta$ equal to $\varphi _j$, $j = 1, 2, \ldots, m$ and substituting them in Eq.~(\ref{sm_weak_form_discretized}), we obtain a system of linear equations of size $m$. We denote it by
\begin{equation}
 A U = F,				\label{matrix_equation_std}
\end{equation}
where $A$ is the $m \times m$ stiffness matrix, $F$ the right-hand side vector of size $m$ and $U$ the solution vector of size $m$, equal to the number of DoFs. 

\subsubsection{Numerical implementation}                               \label{section_solution_technique}
In all the numerical experiments, the IEEE-754 double precision~\cite{zuras2008ieee} is used.
For the time being, we restrict ourselves to two publicly available FEM packages: deal.\rom{2}~\cite{bangerth2007deal} and FEniCS~\cite{alnaes2015fenics}.

\paragraph{Domain discretization}
Unless stated otherwise, we use built-in functions of the FEM packages for generating and refining the computational mesh, and only $h$-regular refinement is considered.
The domain is discretized by different types of elements in deal.\rom{2} and FEniCS: (regular) quadrilaterals are used by the former, and triangles by the latter. 

\begin{figure}[!ht]
  \centering
  \subfloat[deal.\rom{2}\label{sketch_coarsest_computational_mesh_dealii}]{\includegraphics[width=0.15\linewidth]{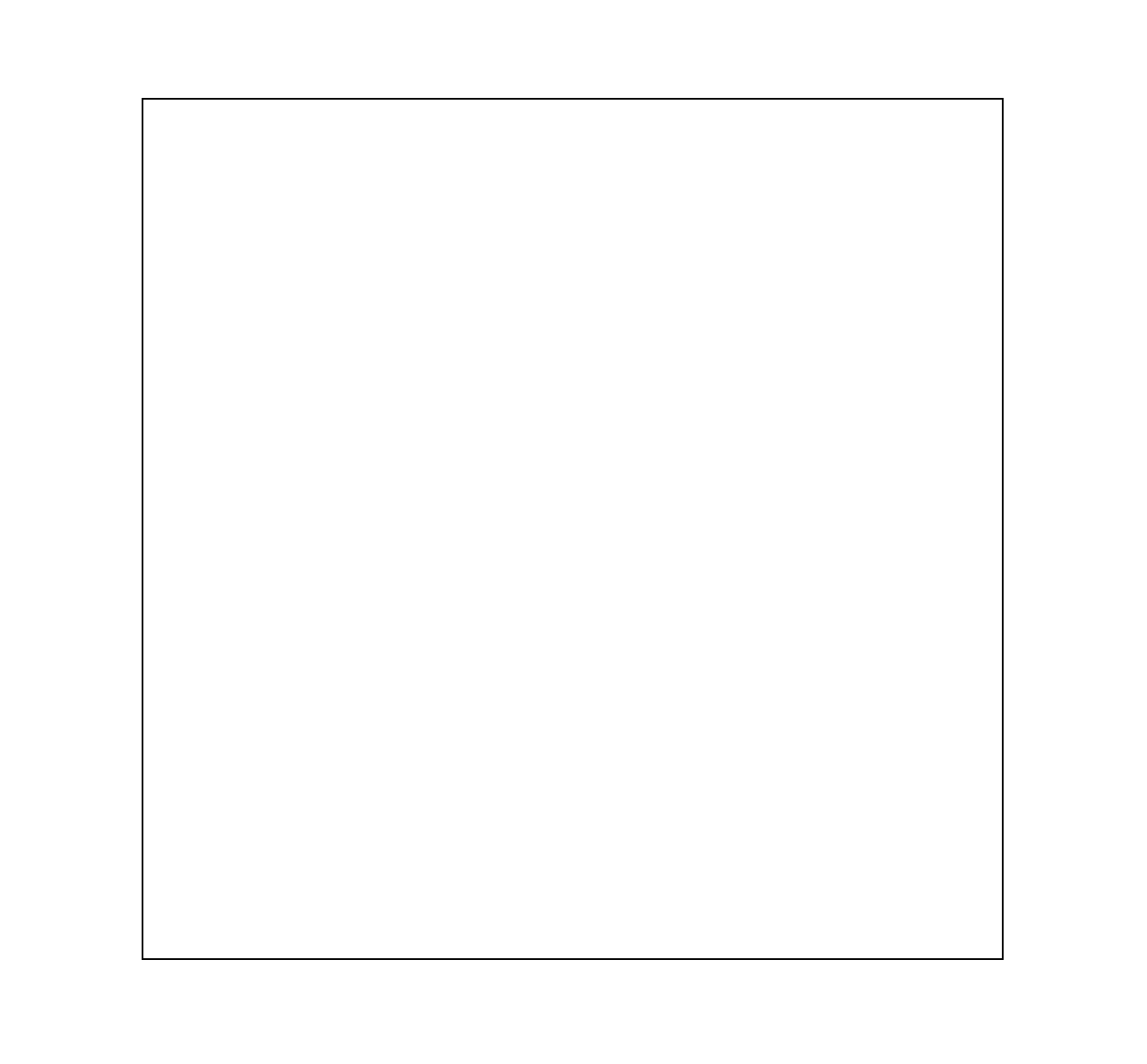}}
  \qquad
  \subfloat[FEniCS\label{sketch_coarsest_computational_mesh_fenics}]{\includegraphics[width=0.15\linewidth]{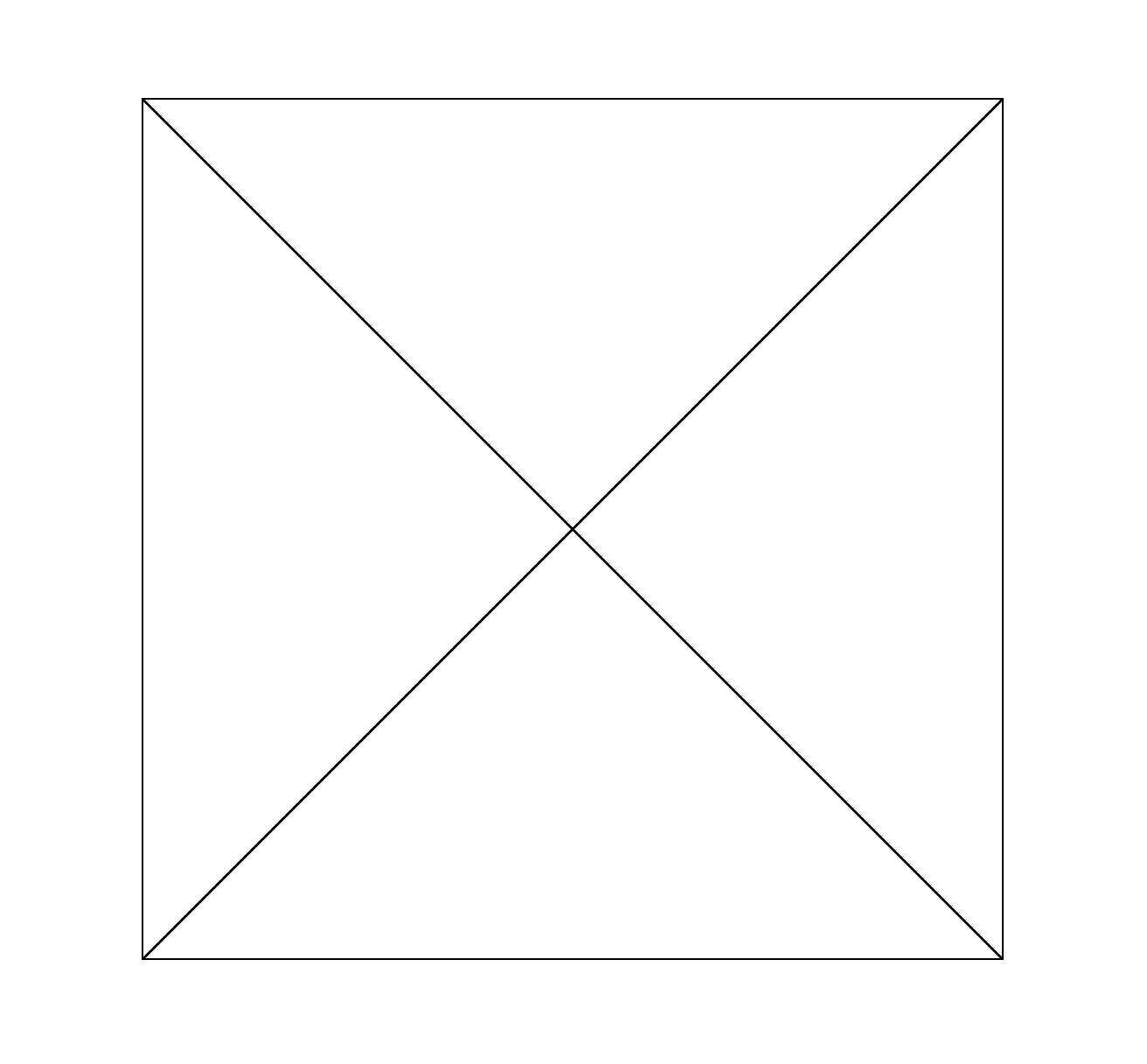}} 
   \caption{Sketch of the original computational mesh.}
   \label{sketch_coarsest_computational_mesh}
\end{figure}

\begin{figure}[!ht]
  \centering
  \subfloat[deal.\rom{2}\label{sketch_computational_mesh_dealii_R_1}]{\includegraphics[width=0.15\linewidth]{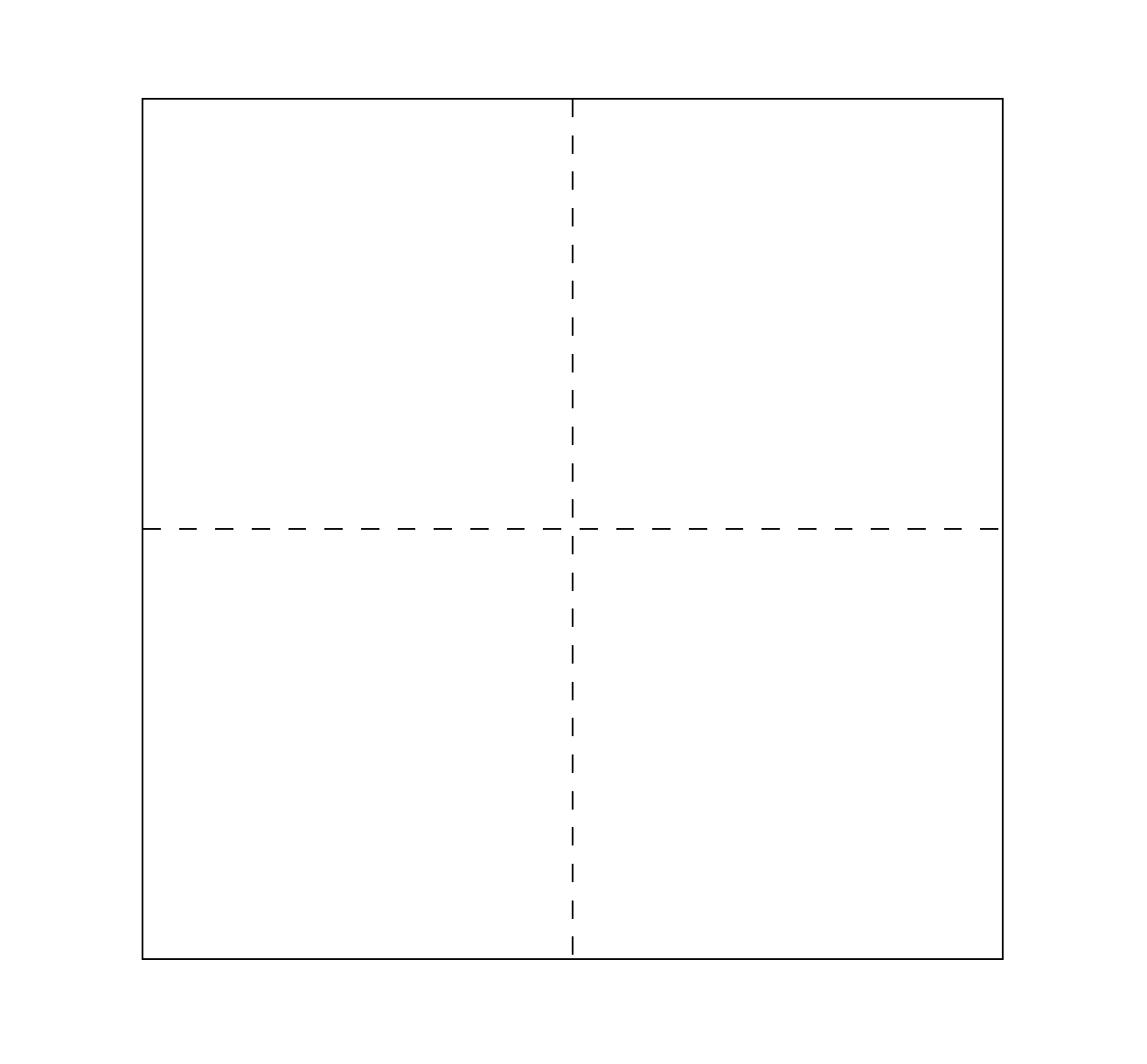}}
  \qquad
  \subfloat[FEniCS\label{sketch_computational_mesh_fenics_R_1}]{\includegraphics[width=0.15\linewidth]{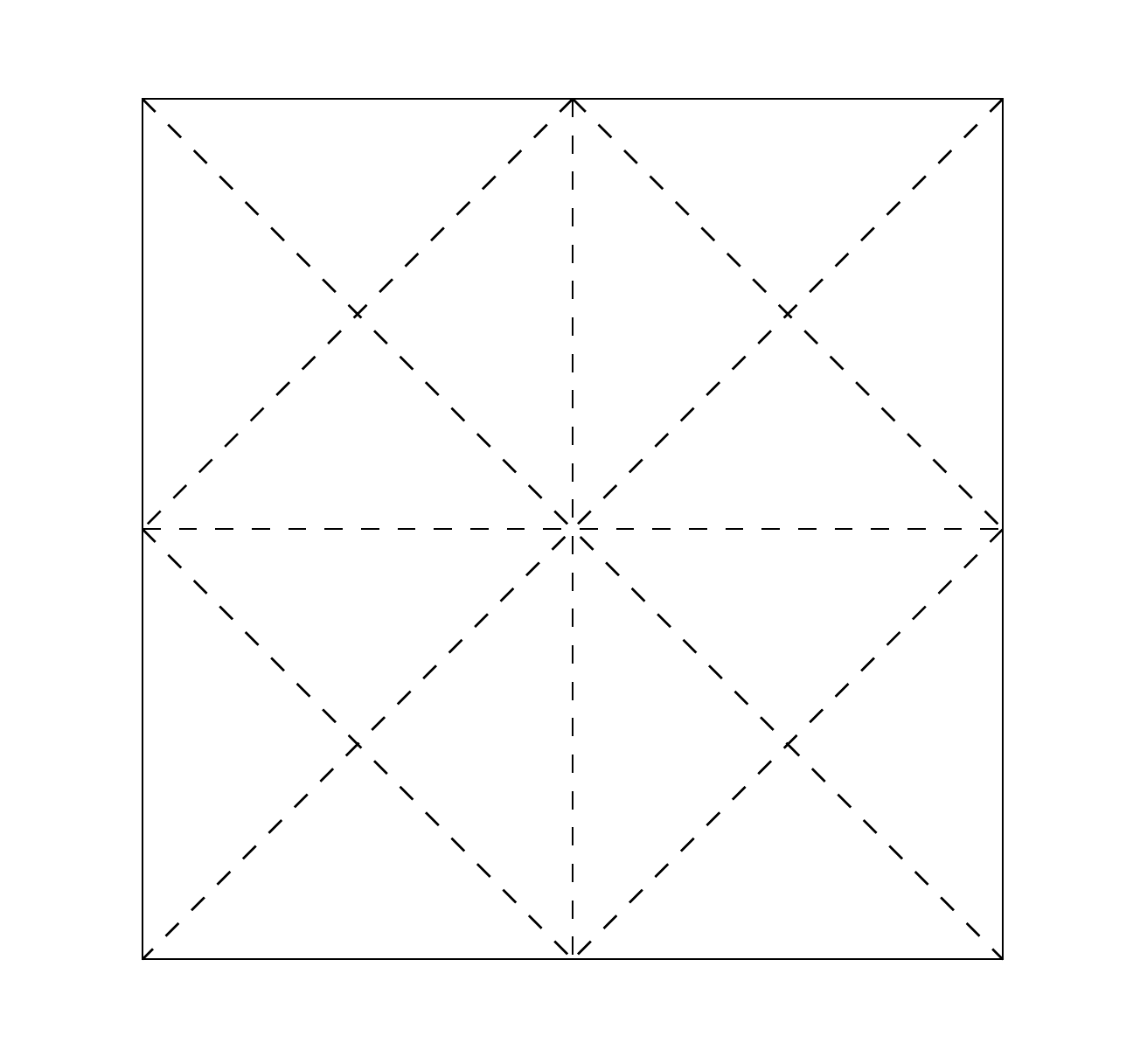}} 
   \caption{Sketch of the computational mesh when $R=1$.}
   \label{sketch_computational_mesh_R_1}
\end{figure}

For deal.\rom{2}, the coarsest computational mesh is shown in Fig.~\ref{sketch_coarsest_computational_mesh_dealii}.
For each $h$-refinement, the mesh is obtained by adding one extra vertex in the center of each element, see Fig.~\ref{sketch_computational_mesh_dealii_R_1} for the computational mesh when the refinement level, denoted by $R$, is 1.
For FEniCS, for each refinement level, the triangle mesh is obtained by further dividing each element of the quadrilateral mesh in deal.\rom{2} into four equal triangles by two diagonals, see Fig.~\ref{sketch_coarsest_computational_mesh_fenics} for the coarsest computational mesh and Fig.~\ref{sketch_computational_mesh_fenics_R_1} for the computational mesh when $R=1$. 
The grid size of a quadrilateral is defined by the side length, and that of a triangle by the height, and hence the grid size using triangles is half of that using quadrilaterals for the same refinement level.

\paragraph{Assembling}
In deal.\rom{2}, the support points are Gauss-Lobatto points, while in FEniCS, the support points are equidistant points, which indicates that the support points of the two packages are different when $p \geq 3$.
For the counting of the number of support points on the computational mesh $K_h$, we refer to \ref{section_counting_number_support_points}.
$T_p$ reads $Q_p$ when using quadrilaterals (as done in deal.\rom{2}) and $P_p$ when using triangles (as done in FEniCS).

\paragraph{Solution method}
To solve the system of equations, the UMFPACK solver~\cite{davis2004algorithm}, which implements the multi-frontal LU factorization approach, is used. Using this solver prevents iteration errors associated with iterative solvers. However, by using this solver, there is an upper limit to the maximum size of the discretized problems because of memory limitations, a problem typically encountered for 2D problems. 
By now, using our hardware, the allowed maximum number of DoFs, which is denoted by $N_{\rm max}$, reads approximately $5 \times 10^7$ and $2 \times 10^6$ for deal.\rom{2} and FEniCS, respectively.

\paragraph{Error estimation}                                                  \label{section_error_estimation}

We investigate the numerical solution $u_h$ obtained with grid size $h$ and its first and second derivatives, see Table~\ref{variables_of_interest_for_computing_the_error}.
Note that for the 2D case, not only the second derivatives $u_{xx}$ and $u_{yy}$ but also the mixed ones, i.e. $u_{xy}$ and $u_{yx}$, come into play.
The error, denoted by $E_h$, is measured in terms of the $L_2$ norm of the difference between the discretized variable and the reference variable.
The reference variable equals the exact expression of the variable, when the exact expression is available or the discretized variable with grid size $h/2$, otherwise~\cite{Runborg2012VerifyingNC}.
We use the exact expression as the reference variable if not stated otherwise.
The solution $u$ is discussed first; next the derivatives are considered.
\begin{table}[!ht]
\caption{Variables of which the error is investigated.}
\label{variables_of_interest_for_computing_the_error}
\centering
\small
\begin{tabular}{l|l|l}      \hline
 Variable & \makecell[c]{1D} & \makecell[c]{2D} \\ \hline
 Solution & $u$ & $u$ \\ \hline
 First derivative & $u_x$ &  $ \nabla u = \begin{bmatrix} u_x \\ u_y \end{bmatrix}$ \\ \hline
 Second derivative & $u_{xx}$ & $\mathbf{H} u = \begin{bmatrix} u_{xx}, u_{xy} \\ u_{yx}, u_{yy} \end{bmatrix}$ \\ \hline
\end{tabular}
\end{table}

\noindent The error of $u_h$ reads
\begin{equation}
  E_h = \| u_h - \widehat{u}\|_2,  \label{relation_error_whole_domain_each_cell}
\end{equation}
where $\| \cdot \|_2$ denotes the $L_2$ norm of a function, and $\widehat{u}$ stands for the reference solution.
The derivatives in 1D, i.e. $u_x$ and $u_{xx}$, only involve one component, and hence we need to replace $u$ in Eq.~(\ref{relation_error_whole_domain_each_cell}) by its derivatives for computing the error;
since the derivatives in 2D, i.e. $\nabla u$ and $\mathbf{H} u$, contain multiple components, we first compute the error of each component according to Eq.~(\ref{relation_error_whole_domain_each_cell}) and then obtain their $l_2$ norm, which denotes the 2-norm of a number sequence.
Note that the error of the second derivative of $u$ does not exist when $p=1$.

\subsubsection{Order of convergence}
Using the error $E_h$ of the discretized variable, the correctness of the implementation in Section \ref{section_solution_technique} is validated by the order of convergence, denoted by $q_h$~\cite{bradji2007convergence}. It is defined by
\begin{equation}
  q_h = \frac{\log \left( \frac{E_{h/2}}{E_{h}} \right)}{\log 2},  \label{formula_order_of_convergence}
\end{equation}
where $E_{h/2}$ is the error of the discretized variable with grid size $h/2$.
According to our experiments and \cite[P.~107]{gockenbach2006understanding}, the asymptotic value of $q_h$, denoted by $q$, reads $p+1$ for the solution and $p + 1 - k$ for the derivatives, where $k$ denotes the order of the derivatives. 
The above $q$ will be used in Section~\ref{section_error_evolution_and_prediction} below for developing the strategy for predicting the highest achievable accuracy.

\subsection{Highest achievable accuracy}                       \label{section_error_evolution_and_prediction}

\subsubsection{Relation between the error and the number of DoFs}                            \label{section_error_evolution_with_N}
\begin{figure}[!ht]
\centering
   \includegraphics[width=0.5\linewidth]{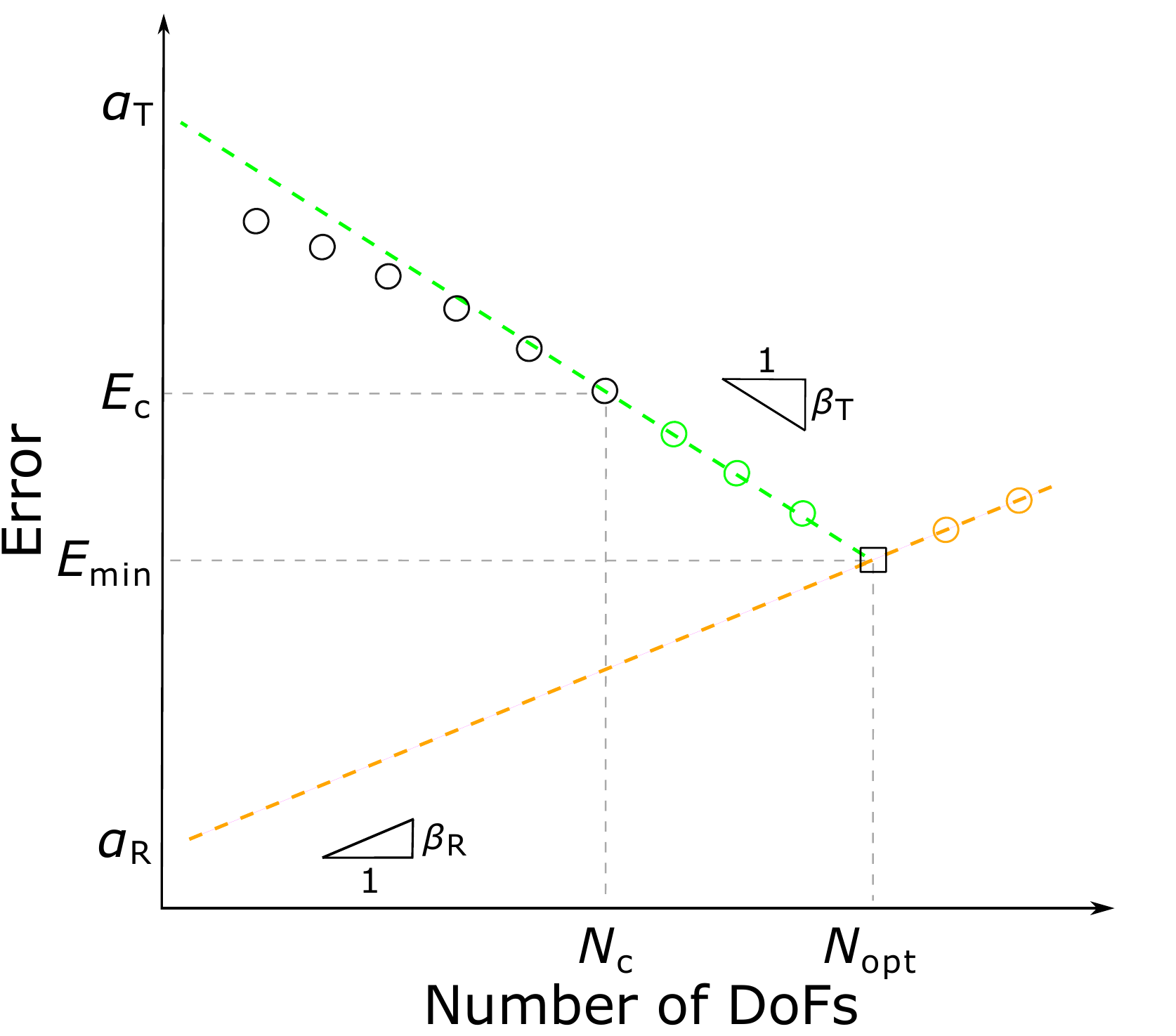}
   \caption{Conceptual sketch of the various errors as a function of the number of DoFs.}
   \label{error_evolution_one_p}
\end{figure}

A conceptual sketch of the error $E_h$ as a function of the number of DoFs, denoted by $N$, in a log-log plot is shown in Fig.~\ref{error_evolution_one_p}, also see \cite{butcher2016numerical}. 
The reason we use the number of DoFs as the $x$ axis, instead of using the grid size $h$, is that the round-off error are independent of $p$~\cite{liu386balancing}.
The error line in Fig.~\ref{error_evolution_one_p} can be divided into a decreasing phase and an increasing one.

In the decreasing phase, the total error $E_h$ is dominated by the truncation error $E_{\rm T}$. 
Based on whether reaching the asymptotic order of convergence $q$ or not, this phase can be further divided into two phases.
When $q$ is not reached (Phase 1), the error path is typically not a straight line, see the black circles.
When $q$ is reached (Phase 2), which often takes a few refinement steps, the error path is a straight line, indicated by the green circles. 
In this phase, $E_h$ can be approximated by
\begin{equation}
 E_{h} \approx E_{\rm T} = \alpha_{\rm T} \cdot {N}^{-\beta_{\rm T}},              \label{formula_truncation_error}
\end{equation}
where $\alpha_{\rm T}$ is the offset of this line, and $\beta_{\rm T}$ the slope, obtained by
\begin{equation*}
  \beta_{\rm T} = - \frac{\log \left( \frac{E_{h/2}}{E_{h}} \right)}{\log \left( \frac{N_{h/2}}{N_{h}} \right)}.
\end{equation*}
In this expression, $N_{h}$ and $N_{h/2}$ are the numbers of DoFs corresponding to $E_h$ and $E_{h/2}$, respectively.
Using Eq.~(\ref{formula_order_of_convergence}), we find
\begin{equation}
  \beta_{\rm T} = \frac{\log 2}{\log \left( \frac{N_{h/2}}{N_{h}} \right)} \cdot q.
\end{equation}

Furthermore, $\frac{N_{h/2}}{N_{h}}$ is about 2 for 1D problems and 4 for 2D problems. 
The proof for the latter can be found in \ref{section_proof_slope_ET}, and that for the former follows in the similar manner. 
Therefore, we have $\beta_{\rm T} \approx q$ for 1D problems and $\beta_{\rm T} \approx q/2$ for 2D problems.

Denoting the number of DoFs and the error when $q$ is reached as $N_c$ and $E_{\rm c}$, respectively, substituting them into Eq.~(\ref{formula_truncation_error}), we have
\begin{equation}
 \alpha_{\rm T} = {E_{\rm c}} \cdot {N_{\rm c}}^{\beta_{\rm T}}.		\label{formula_offset_truncation_error}
\end{equation}

In the increasing phase, the total error $E_h$ is dominated by the round-off error $E_{\rm R}$.
As suggested in \cite{liu386balancing,Babuska2018Roundoff}, the error in this phase tends to increase constantly, see the orange circles in Fig.~\ref{error_evolution_one_p}. Denoting the offset by $\alpha_{\rm R}$ and the slope by $\beta _{\rm R}$, one finds that
\begin{equation}
 E_{h} \approx E_{\rm R} = \alpha_{\rm R} \cdot {N}^{\beta_{\rm R}},                  \label{formula_round_off_error}
\end{equation}
with the coefficients $\alpha_{\rm R}$ and $\beta _{\rm R}$ to be determined in Section \ref{section_novel_approach_for_obtaining_round_off_error}.
At the turning point of the decreasing phase and the increasing phase, denoted by the black square, the minimum error $E_{\rm min}$ is obtained.
The behaviour of $E_h$ as a function of the number of DoFs is summarized in Table~\ref{table_error_in_different_phases}, and depicted in Fig. \ref{error_evolution_one_p}.
\begin{table}[!ht]
\small
\caption{Description of the evolution of $E_h$.}
\label{table_error_in_different_phases}
\centering
\scriptsize
 \begin{tabular}{l | c | c | c} \hline
 & \multicolumn{2}{c|}{Decreasing phase} & \multirow{2}{*}{Increasing phase} \\ \cline{2-3} 
 & Phase 1 & Phase 2 & \\ \hline
Size of $N$  & $N < N_{\rm c}$ & $N_{\rm c} \leqslant N < N_{\rm opt}$ & $N_{\rm opt} \leqslant N$ \\ \hline
Phenomenon & \makecell[l]{Decreasing but not\\ converging with slope $\beta_{\rm T}$} & \makecell[l]{Decreasing and converging\\ with slope $\beta_{\rm T}$ and offset $\alpha_{\rm T}$} & \makecell[l]{Increasing and converging\\ with slope $\beta_{\rm R}$ and offset $\alpha_{\rm R}$} \\	\hline
Dominant error & \multicolumn{2}{c|}{Truncation error} & Round-off error \\	\hline
Formula & - & $E_h \approx E_{\rm T}=\alpha_{\rm T}{N}^{-\beta_{\rm T}}$ & $E_h \approx E_{\rm R}=\alpha_{\rm R} {N}^{\beta _{\rm R}}$ \\	\hline
\end{tabular}
\end{table}

\subsubsection{Obtaining the highest achievable accuracy}                               \label{section_scheme_for_predicting_the_minimum_error}
There are two methods for obtaining the highest achievable accuracy.
One is the brute-force method, denoted by BF.
This method uses a sequential number of numerical experiments until the error starts to increase (indicated by the black and green circles, black square, and first orange circle in Fig.~\ref{error_evolution_one_p}), which is time-consuming and impractical. 

\begin{subequations}
The other is the method suggested in \cite{liu386balancing}, which is denoted by PRED+.
In this method, first, when the asymptotic order of convergence $q$ is reached, $N_{\rm opt}$ and $E_{\rm min}$ are predicted as follows:
\begin{align}
 N_{\rm opt} &= \left( \frac{\alpha_{\rm T} \cdot \beta_{\rm T}}{\alpha _{\rm R} \cdot \beta_{\rm R}} \right)^{\frac{1}{\beta_{\rm T} + \beta_{\rm R}}}, \label{formula_N_opt}      \\
 E_{\rm min} &= \alpha_{\rm T} \cdot {N_{\rm opt}}^{- {\beta _{\rm T}}}+\alpha_{\rm R} \cdot {N_{\rm opt}}^{{\beta _{\rm R}}},        \label{formula_E_min} 
\end{align}
\end{subequations}
Next, the solution with the highest achievable accuracy is obtained by computing the result using the $N_{\rm opt}$ predicted.
Since a number of $h$-refinements are circumvented using the PRED+ method, the CPU time required for obtaining the discretized variable with the accuracy $E_{\rm min}$, denoted by $T$, using the  PRED+ method is expected to be much less than that using the BF method.

\section{Novel method for determining \texorpdfstring{$\alpha_{\rm R}$}{alphaR} and \texorpdfstring{$\beta_{\rm R}$}{betaR} of the round-off error}                \label{section_novel_approach_for_obtaining_round_off_error}        

\subsection{Strategy}

In \cite{liu386balancing}, $\alpha_{\rm R}$ and $\beta_{\rm R}$ are fitted from the orange circles illustrated in Fig.~\ref{error_evolution_one_p}.
It showed that for many specific 1D problems, using the standard FEM in deal.\rom{2}, when the $L_2$ norm of the solution, i.e. $\| u \|_2$, is of order 1, the offset $\alpha_{\rm R}$ is about $2 \times 10^{-17}$ for $u$ and increases slightly with the increasing order of derivative;
the slope $\beta_{\rm R}$ is constant, equals 2, and is independent of the variable. 
Furthermore, for different model problems, $\alpha_{\rm R}$ is linearly proportional to $\| u \|_2$.
Using the above information, the coefficients $\alpha_{\rm R}$ and $\beta_{\rm R}$ were estimated for general 1D problems when solving them using the standard FEM in deal.\rom{2}.
Using the above coefficients $\alpha_{\rm R}$ and $\beta_{\rm R}$ in the PRED+ method, the accuracy $E_{\min}$ predicted is very close to that obtained using the BF method. 

However, $\alpha_{\rm R}$ and $\beta_{\rm R}$ might be different for higher space dimensions, different types and packages of FEM.
Even though fewer problem cases can be chosen to obtain the orange circles, the CPU time required is still very large since the numbers of DoFs corresponding to the orange circles are larger than $N_{\rm opt}$.

\begin{figure}[!ht]
\centering
   \includegraphics[width=0.5\linewidth]{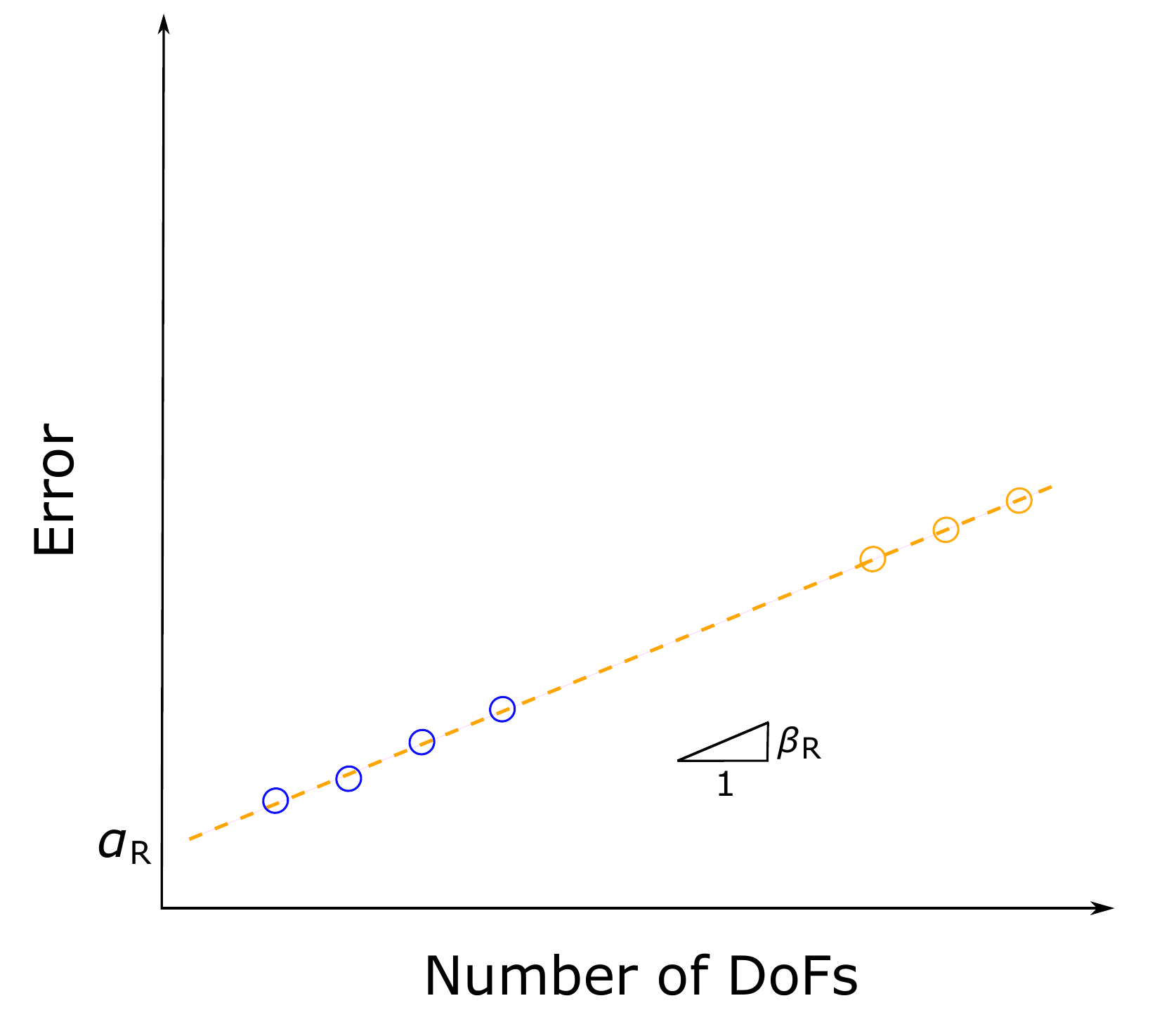}
   \caption{Illustration of using the MS+ method for obtaining $\alpha_{\rm R}$ and $\beta_{\rm R}$.}
   \label{sketch_illustration_mms_one_p}
\end{figure}

Fortunately, the round-off error can also be obtained with a few cheap experiments when the solution of a problem can be exactly represented in the FE space under consideration (indicated by the blue circles in Fig.~\ref{sketch_illustration_mms_one_p})~\cite{liu386balancing}. 
Therefore, using the same settings as for the problem at hand, i.e. $\Omega$, $D(\mathbf{x})$, $r(\mathbf{x})$, and type of boundary conditions, we propose to first obtain $\alpha_{\rm R}$ and $\beta_{\rm R}$ for a problem with a manufactured solution $u_{\rm M}$ that can be exactly represented in the current FE space.
$\alpha_{\rm R}$ and $\beta_{\rm R}$ for $u_{\rm M}$ are denoted by $\alpha_{\rm R, M}$ and $\beta_{\rm R, M}$, respectively.

Next, since $\| u_{\rm M} \|_2$ does not have to equal $\| u_{\rm O} \|_2$, where $u_{\rm O}$ denotes the original solution, we have to adjust $\alpha_{\rm R, M}$ according to the linear relation between $\alpha_{\rm R}$ and $\| u \|_2$.
$\alpha_{\rm R}$ adjusted from $\alpha_{\rm R, M}$, denoted by $\alpha_{\rm R, M+}$, reads
\begin{equation}
  \alpha_{\rm R, M+} = \alpha_{\rm R, M} \cdot \frac{\| u_{\rm O} \|_2}{\| u_{\rm M} \|_2}.           \label{formula_alpha_R_prediction}
\end{equation}

We call this method the method of manufactured solutions assisted by $\| u \|_2$, denoted by MS+~\cite{salari2000code}, and the method that uses the orange circles the method of original solutions,       denoted by OS. Using the OS method, the resulting $\alpha_{\rm R}$ and $\beta_{\rm R}$ are denoted by $\alpha_{\rm R, O}$ and $\beta_{\rm R, O}$, respectively.
In Section~\ref{section_validation_MS_plus} below, we will test the MS+ method by comparing $\alpha_{\rm R, M+}$ and $\beta_{\rm R, M}$ with $\alpha_{\rm R, O}$ and $\beta_{\rm R, O}$.

\subsection{Results}				\label{section_validation_MS_plus}

In this section, we investigate the accuracy of approximating $\alpha_{\rm R, O}$ and $\beta_{\rm R, O}$ by $\alpha_{\rm R, M+}$ and $\beta_{\rm R,M}$, respectively, by comparing their values.
All problems will be solved using the deal.\rom{2} software, and the element degree ranges from 1 to 5.
In Section~\ref{section_1d_benchmark_problems}, 1D problems are considered, followed by 2D problems in Section~\ref{section_2d_benchmark_problems}.

\subsubsection{1D problems}                    \label{section_1d_benchmark_problems}

For the 1D case, we investigate problems with $u=e^{-(x-0.5)^2}$, which cannot be reproduced exactly in the FE space $V_g$.
Appropriate Dirichlet boundary conditions are imposed on both ends.
The manufactured solution $u_{\rm M}$ is chosen to be $(x-0.5)^2$.
$\|u_{\rm O}\|_2 \approx 0.92$ is about 8.4 times $\|u_{\rm M}\|_2 \approx 0.11$, and hence $\alpha_{\rm R, M+} = 8.4 \cdot \alpha_{\rm R,M}$ using Eq.~(\ref{formula_alpha_R_prediction}).

\paragraph{Poisson problems}			\label{section_1d_benchmark_problems_pois}

We first investigate results obtained using an equidistant mesh and then consider the influence of distorted meshes.

\begin{figure}[!ht]
	\subfloat[$u$\label{error_evolution_1d_0_pois_0p0_u_x_m_0p5_square_uniform_0_sm_solu}]{\includegraphics[width=0.33\linewidth]{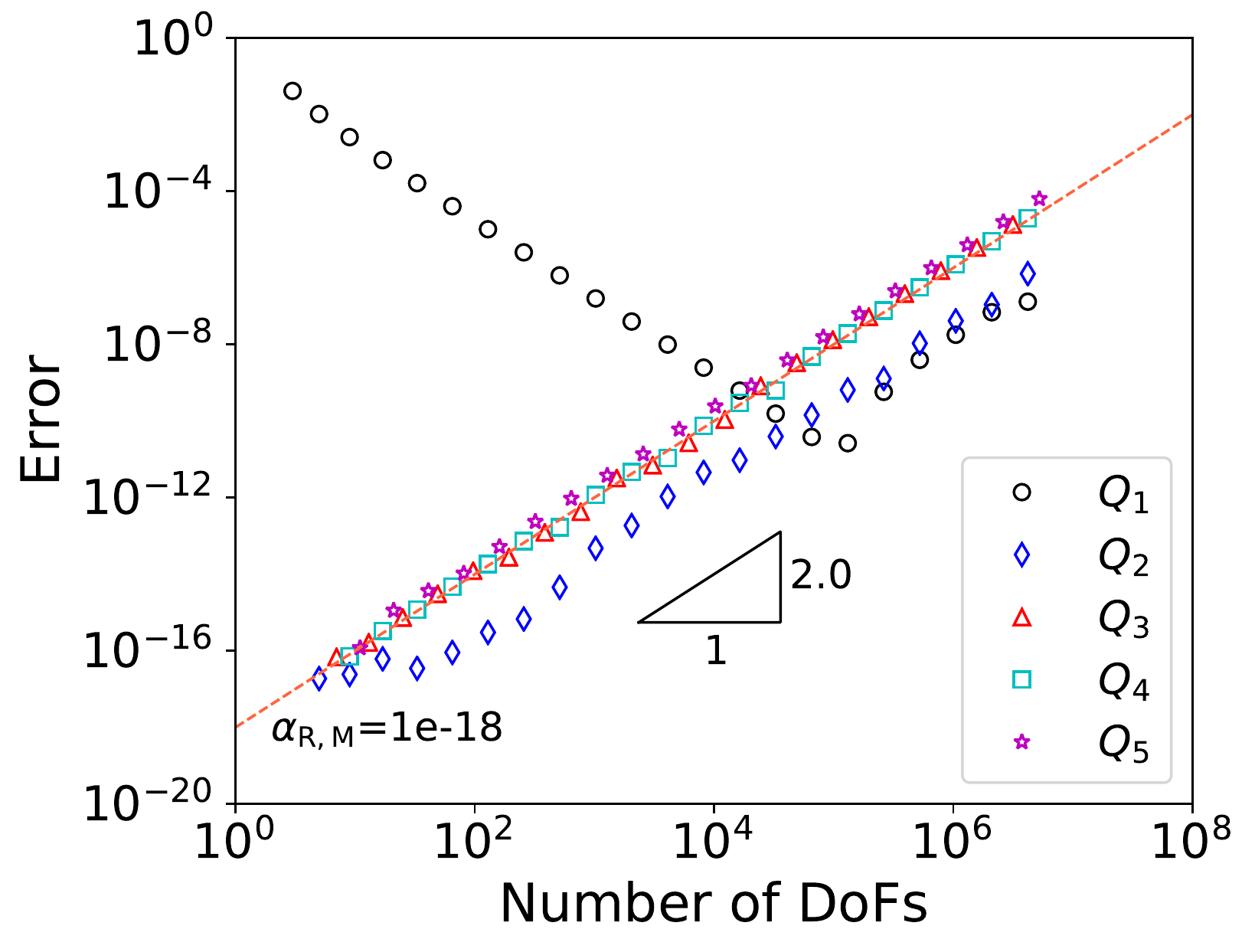}}
    \subfloat[$u_x$\label{error_evolution_1d_0_pois_0p0_u_x_m_0p5_square_uniform_0_sm_grad}]{\includegraphics[width=0.33\linewidth]{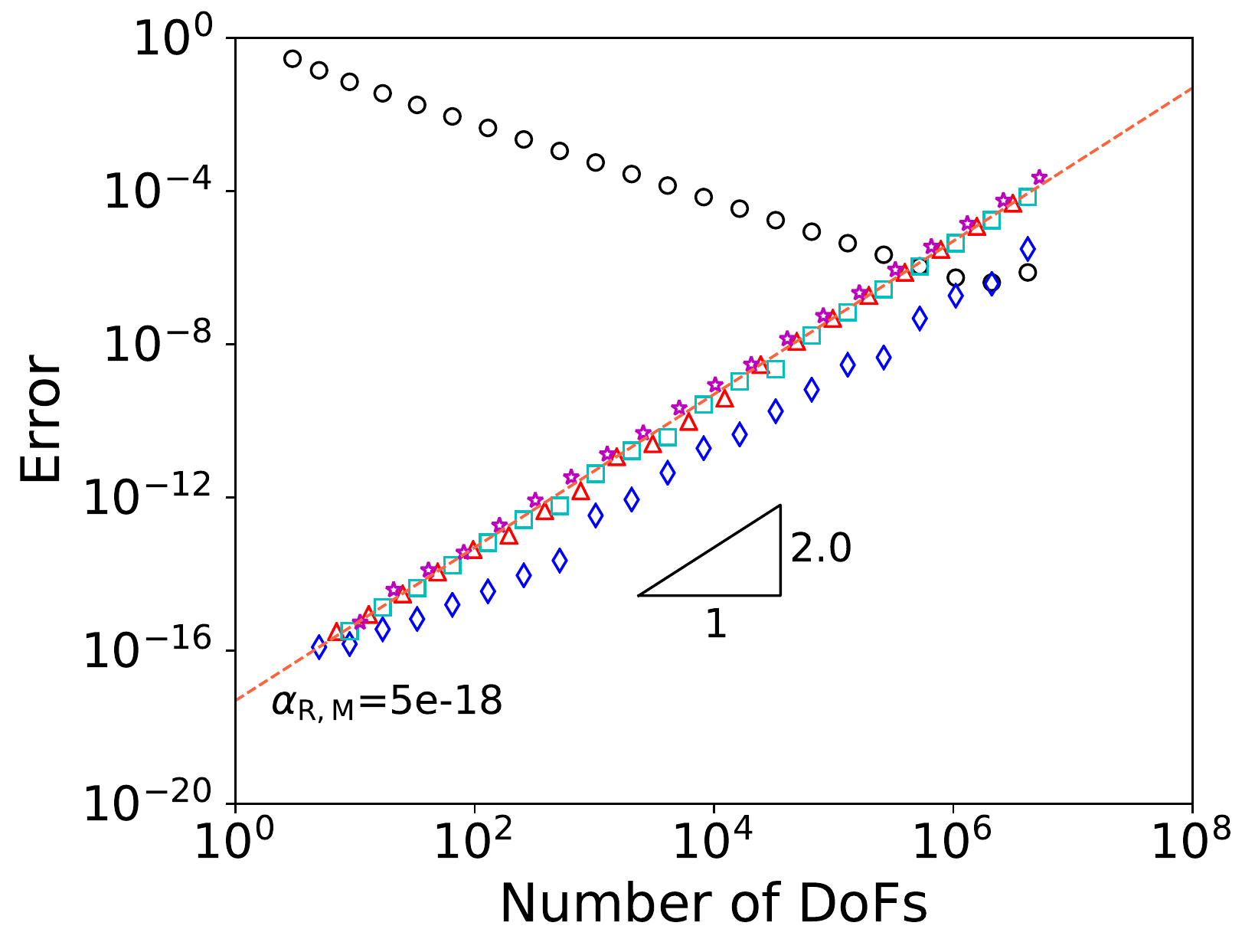}}
    \subfloat[$u_{xx}$\label{error_evolution_1d_0_pois_0p0_u_x_m_0p5_square_uniform_0_sm_2ndd}]{\includegraphics[width=0.33\linewidth]{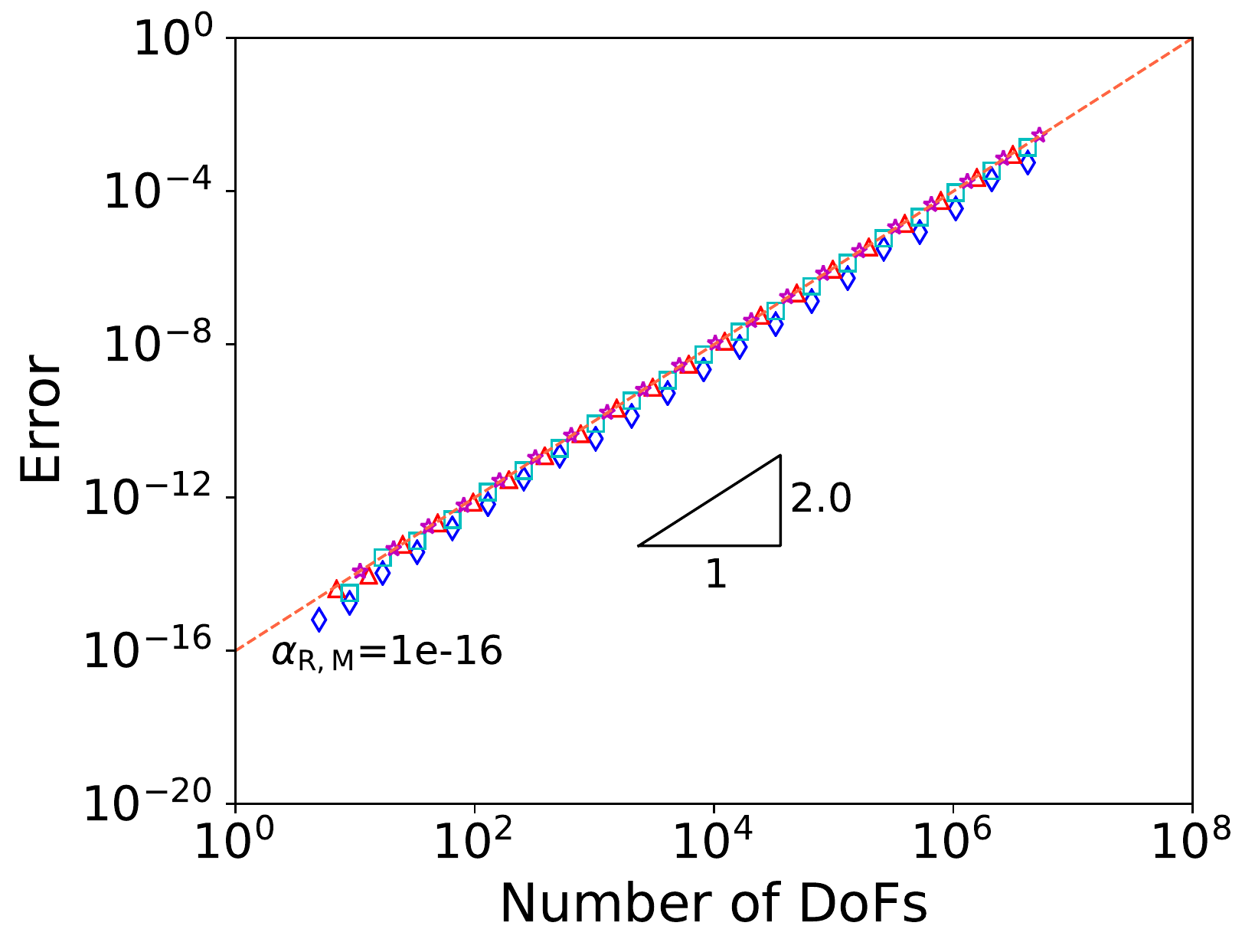}}
\caption{Error evolution for the 1D Poisson problem with $u = (x-0.5)^2$ and Mesh Type 1.}
\label{error_evolution_1d_0_pois_0p0_u_x_m_0p5_square_uniform_0_sm}
\end{figure}

Using an equidistant mesh, the dependency of the round-off error for $u_{\rm M}$ is shown in Fig.~\ref{error_evolution_1d_0_pois_0p0_u_x_m_0p5_square_uniform_0_sm}.
From this figure, it follows that as expected for $p > 1$, the only error source is the round-off error.
For all variables $u$, $u_x$ and $u_{xx}$, the round-off error line is a straight line in the log-log plot, with slope $\beta_{\rm R,M} = 2$, while the offset $\alpha_{\rm R,M}$ depends on the variable under consideration, see Fig.~\ref{error_evolution_1d_0_pois_0p0_u_x_m_0p5_square_uniform_0_sm} and columns 3 and 4 in the second row of Table~\ref{table_results_using_MS_plus_for_1d_pois_u_exp_m_x_m_0p5_square_dealii} for their specific values.
Using this information, the values of $\alpha_{\rm R, M+}$ obtained are shown in column 5 in the second row of Table~\ref{table_results_using_MS_plus_for_1d_pois_u_exp_m_x_m_0p5_square_dealii}.

\begin{figure}[!ht]
	\subfloat[$u$\label{error_evolution_1d_0_pois_1p0_u_exp_m_x_m_0p5_square_uniform_0_sm_solu}]{\includegraphics[width=0.33\linewidth]{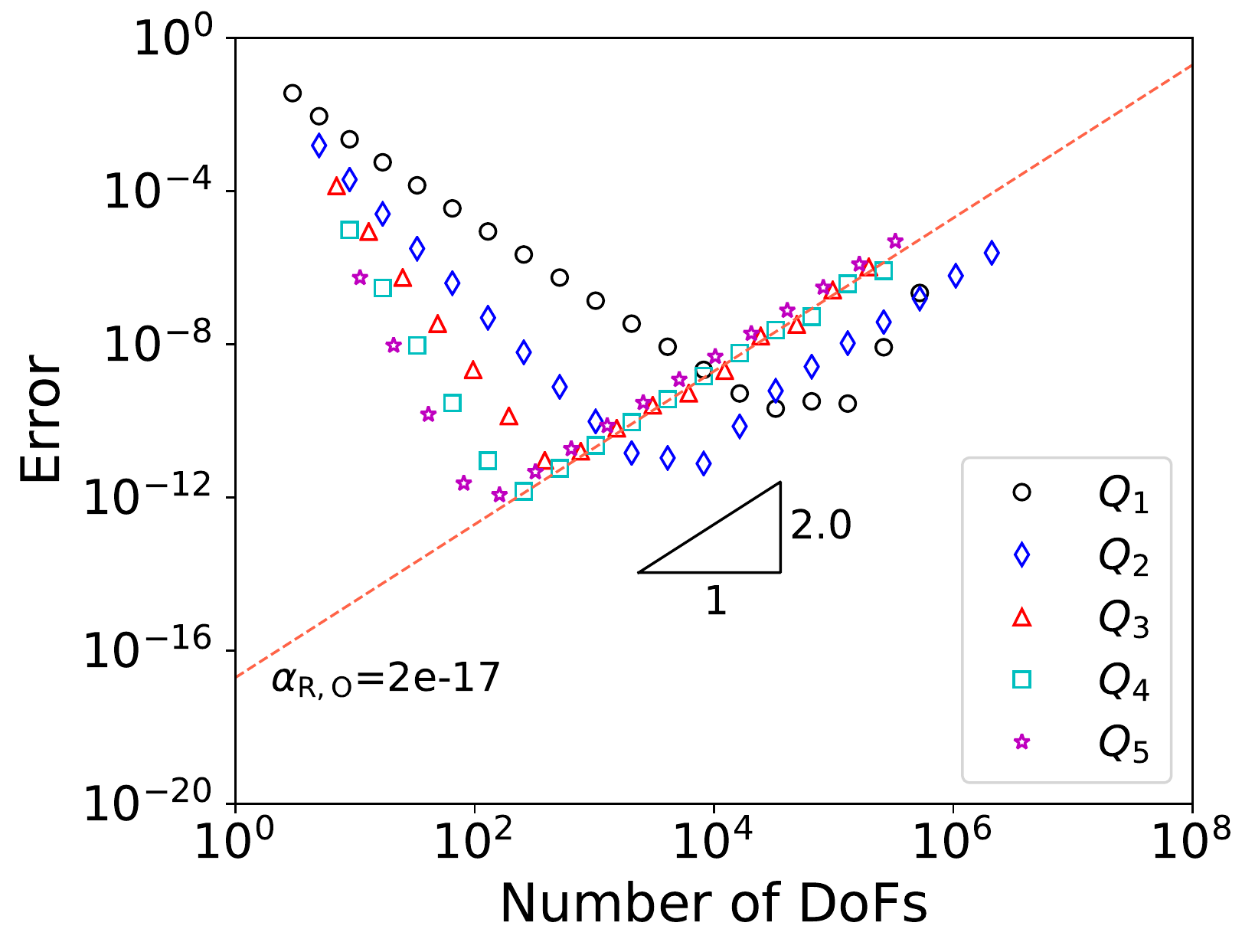}}
    \subfloat[$u_x$\label{error_evolution_1d_0_pois_1p0_u_exp_m_x_m_0p5_square_uniform_0_sm_grad}]{\includegraphics[width=0.33\linewidth]{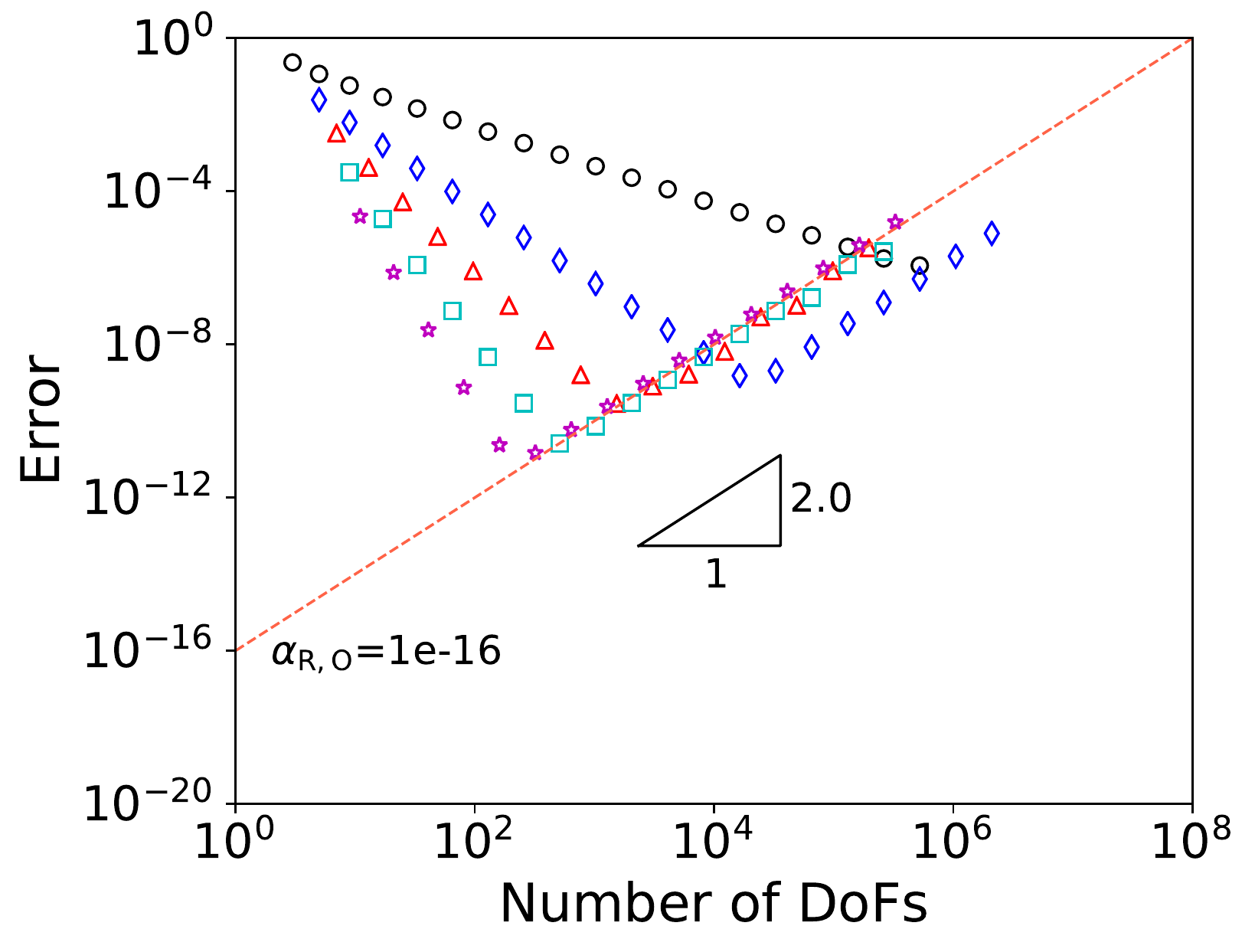}}
    \subfloat[$u_{xx}$\label{error_evolution_1d_0_pois_1p0_u_exp_m_x_m_0p5_square_uniform_0_sm_2ndd}]{\includegraphics[width=0.33\linewidth]{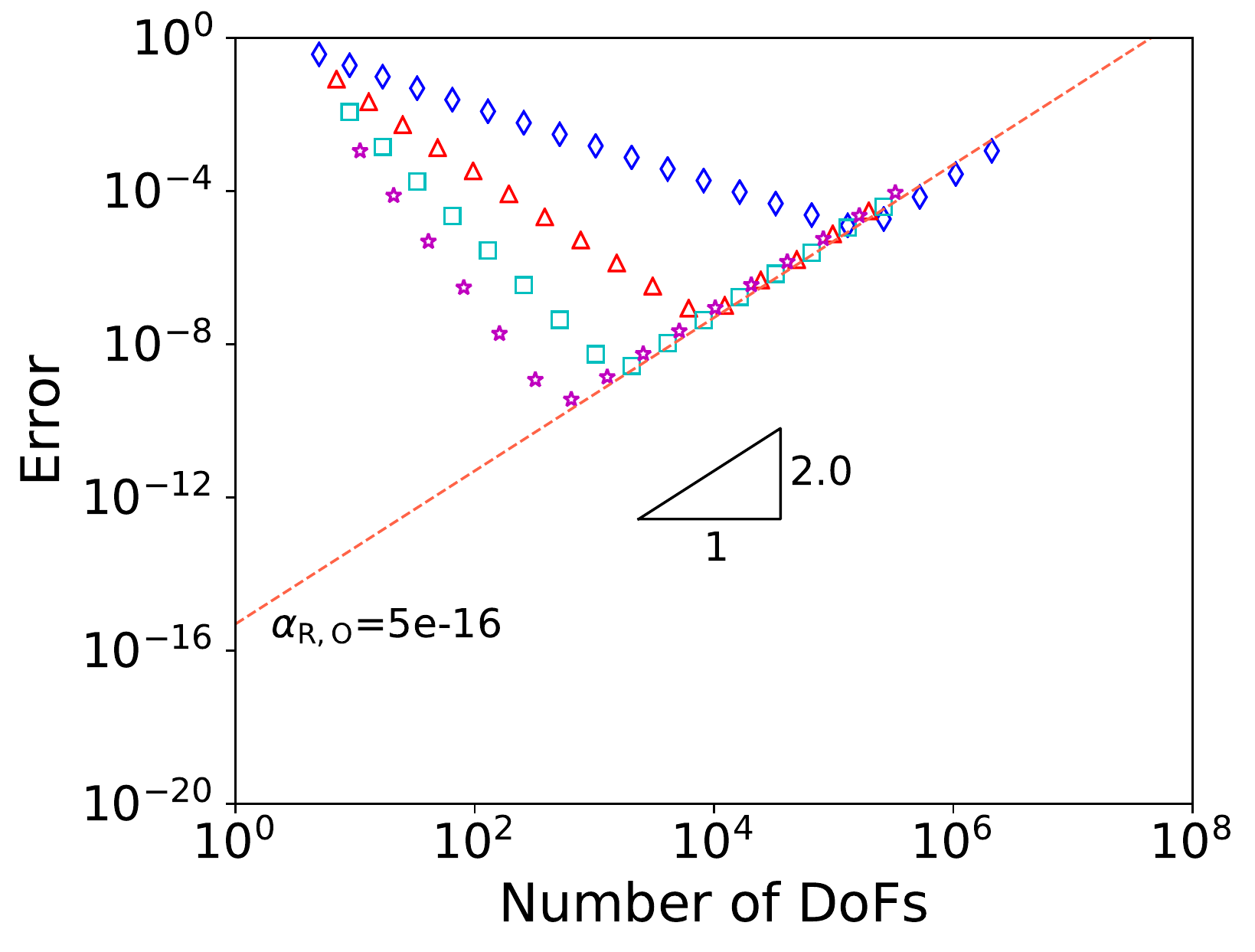}}
\caption{Error evolution for the 1D Poisson problem with $u=e^{-(x-0.5)^2}$ and Mesh Type 1.}
\label{error_evolution_1d_0_pois_1p0_u_exp_m_x_m_0p5_square_uniform_0_sm}
\end{figure}

Using the OS method, the dependency of the error is shown in Fig.~\ref{error_evolution_1d_0_pois_1p0_u_exp_m_x_m_0p5_square_uniform_0_sm}.
As can be seen, the truncation error tends to converge at the analytical order for different $p$ before reaching the round-off error.
Analyzing the round-off error line, it follows that $\beta_{\rm R, O} = 2$, and $\alpha_{\rm R, O}$ increases when increasing the order of the derivative, see Fig.~\ref{error_evolution_1d_0_pois_1p0_u_exp_m_x_m_0p5_square_uniform_0_sm} and columns 6 and 7 in the second row of Table~\ref{table_results_using_MS_plus_for_1d_pois_u_exp_m_x_m_0p5_square_dealii} for the specific values.

\begin{table}[!ht]
\caption{The coefficients $\alpha_{\rm R}$ and $\beta_{\rm R}$ obtained using the method of MS+ and OS for the 1D Poisson problem with $u = e^{-(x-0.5)^2}$.}
\centering
\scriptsize
\begin{tabular}{c|c|c|c|c|c|c}	\hline
 \multirow{2}{*}{Mesh type} & \multirow{2}{*}{Variable} & \multicolumn{3}{c|}{MS+} & \multicolumn{2}{c}{OS} \\ \cline{3-7}
   &  & $\alpha_{\rm R,M}$ & $\beta_{\rm R,M}$ & $\alpha_{\rm R, M+}$ & $\alpha_{\rm R, O}$ & $\beta_{\rm R, O}$ \\ \hline
 1 & \makecell{$u$\\ $u_x$\\ $u_{xx}$} & \makecell{1.0e-18\\ 5.0e-18\\ 1.0e-16} & \multirow{7}{*}{\makecell{2.0\\ 2.0\\ 2.0}} & \makecell{8.4e-18\\ 4.2e-17\\ 8.4e-16} & \makecell{2.0e-17\\ 1.0e-16\\ 5.0e-16} & \multirow{10}{*}{\makecell{2.0\\ 2.0\\ 2.0}} \\ \cline{1-3} \cline{5-6}
 2 & \makecell{$u$\\ $u_x$\\ $u_{xx}$} & \multirow{4}{*}{\makecell{1.0e-18\\ 5.0e-18\\ 2.0e-16}} &  & \multirow{4}{*}{\makecell{8.4e-18\\ 4.2e-17\\ 1.7e-15}} & \makecell{2.0e-17\\ 1.0e-16\\ 5.0e-15} &  \\ \cline{1-2} \cline{6-6}
 3 & \makecell{$u$\\ $u_x$\\ $u_{xx}$} &  &  &  & \makecell{2.0e-17\\ 1.0e-16\\ 2.0e-15} &  \\ \cline{1-6}
 4 & \makecell{$u$\\ $u_x$\\ $u_{xx}$} & \makecell{1.0e-18\\ 5.0e-18\\ 1.0e-16} & \makecell{1.8\\ 1.8\\ 2.0} & \makecell{8.4e-18\\ 4.2e-17\\ 8.4e-16} & \makecell{2.0e-18\\ 5.0e-18\\ 1.0e-15} &  \\ \hline
\end{tabular}
\label{table_results_using_MS_plus_for_1d_pois_u_exp_m_x_m_0p5_square_dealii}
\end{table}

\begin{figure}[!ht]
  \centering
	\subfloat[Type 1\label{manufacturing_alpha_R_beta_R_0_pois_Type_1}]{\includegraphics[width=0.33\linewidth]{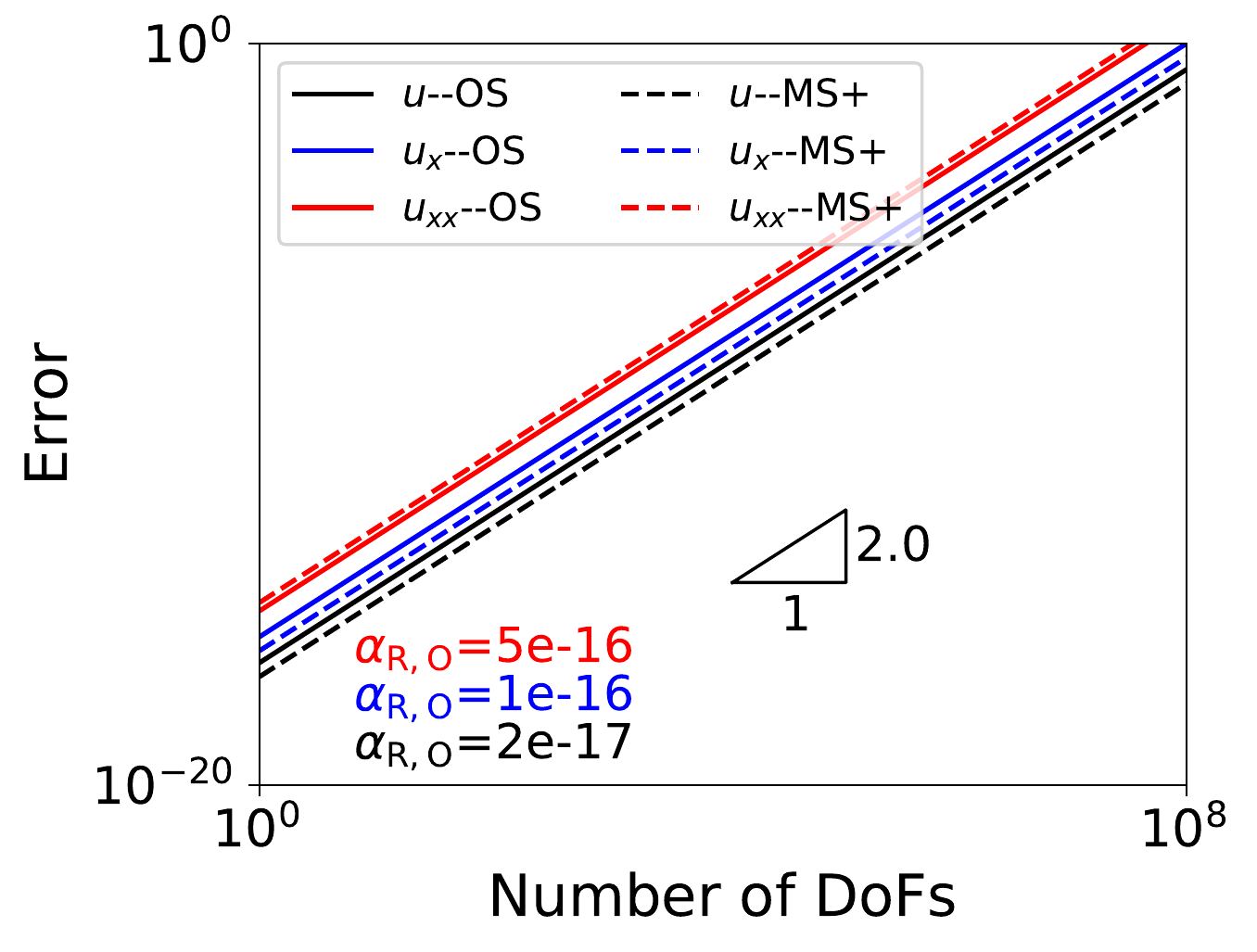}}
    \subfloat[Type 2\label{manufacturing_alpha_R_beta_R_0_pois_Type_2}]{\includegraphics[width=0.33\linewidth]{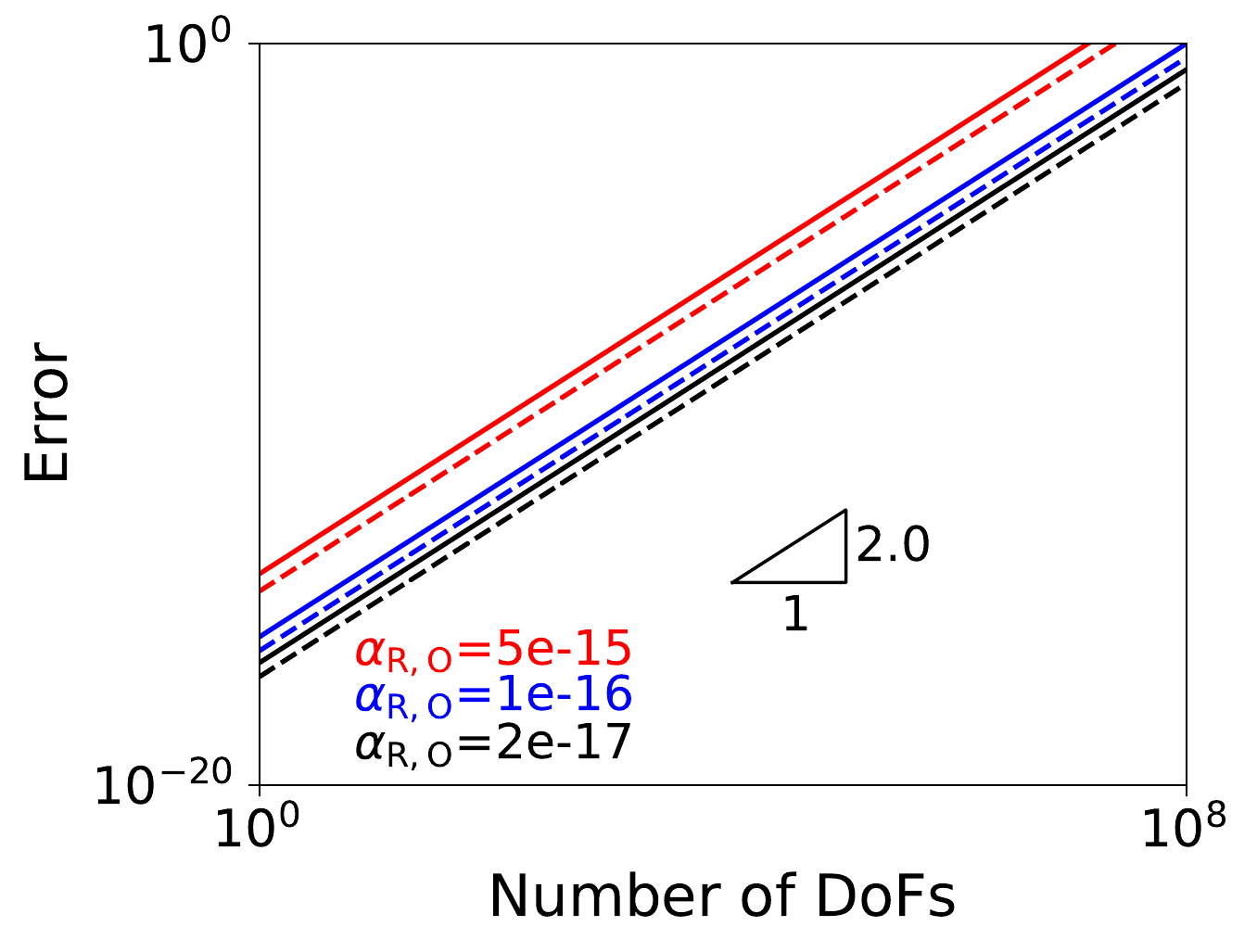}}
    \\
    \subfloat[Type 3\label{manufacturing_alpha_R_beta_R_0_pois_Type_3}]{\includegraphics[width=0.33\linewidth]{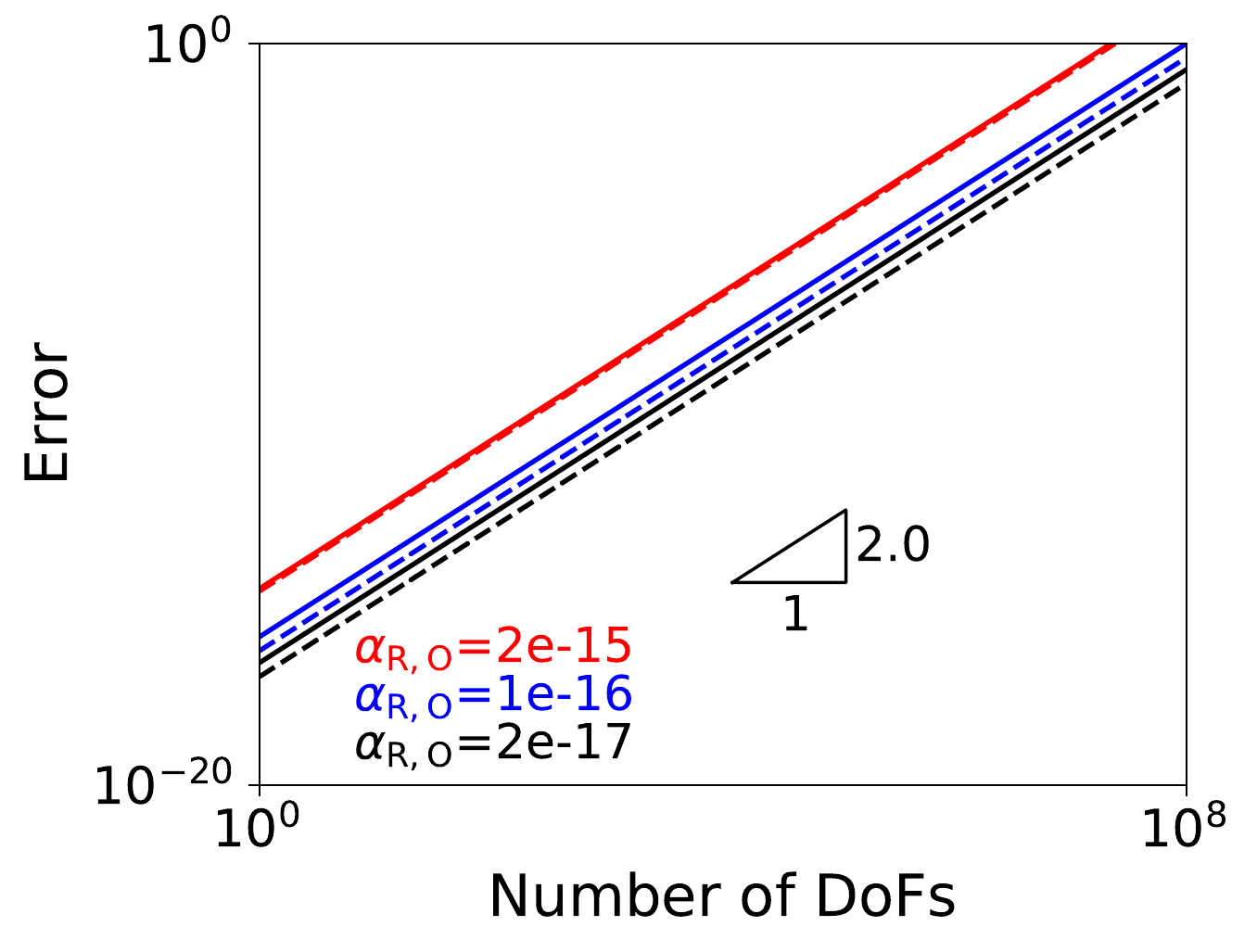}}
    \subfloat[Type 4\label{manufacturing_alpha_R_beta_R_0_pois_Type_4}]{\includegraphics[width=0.33\linewidth]{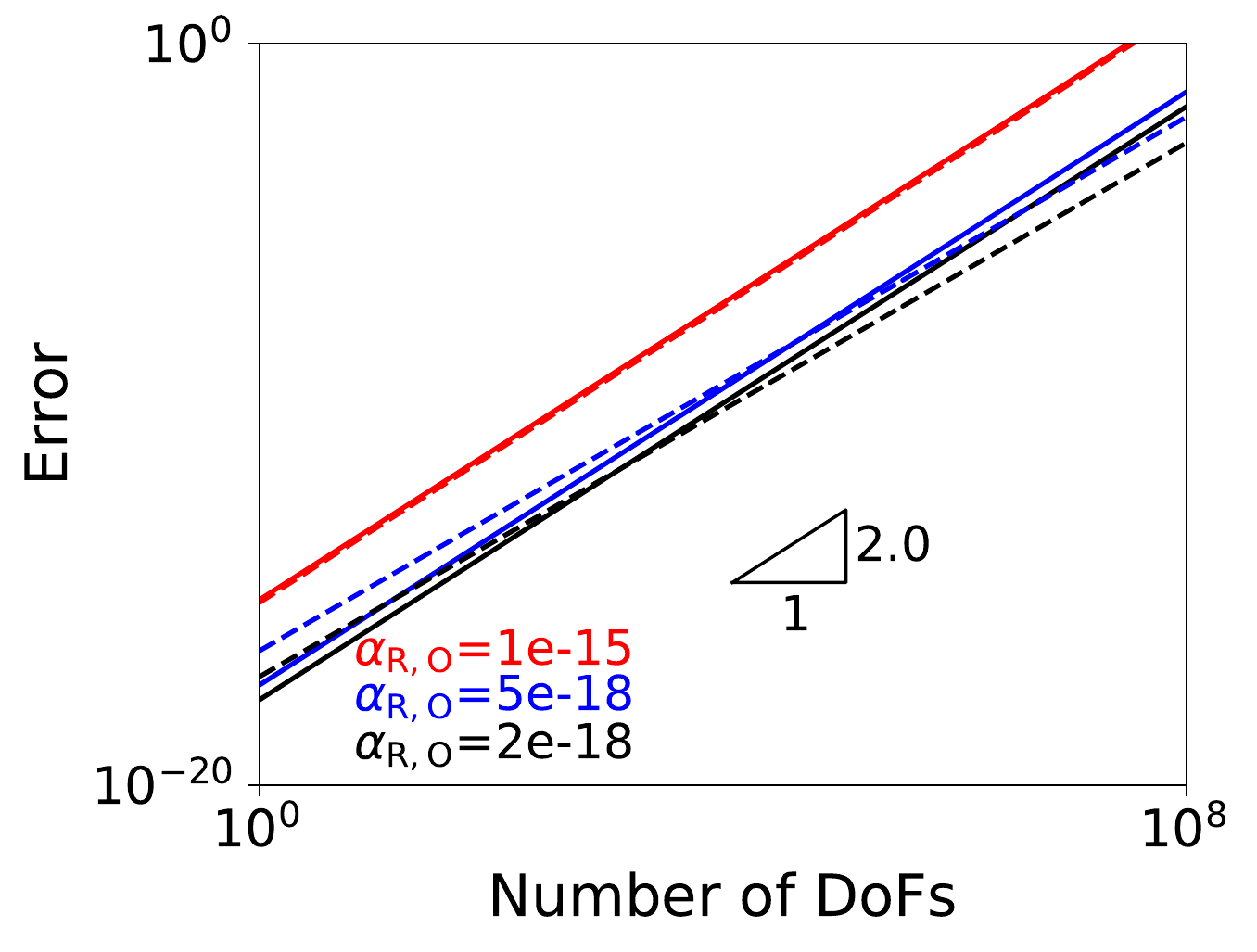}}    
\caption{Comparison of the round-off error line represented by $\alpha_{\rm R, M+}$ and $\beta_{\rm R, M}$ to that represented by $\alpha_{\rm R, O}$ and $\beta_{\rm R, O}$ for the 1D Poisson problem with $u=e^{-(x-0.5)^2}$.}
\label{manufacturing_alpha_R_beta_R_0_pois}
\end{figure}

\newpage
Finally, the round-off error line represented by $\alpha_{\rm R, M+}$ and $\beta_{\rm R, M}$ is compared to that represented by $\alpha_{\rm R, O}$ and $\beta_{\rm R, O}$ in Fig.~\ref{manufacturing_alpha_R_beta_R_0_pois_Type_1}. 
In this figure, the solid lines correspond to the round-off error using the OS method, while the dashed lines correspond to the round-off error using the MS+ method;
the results for the principle variable $u$ and derivatives $u_x$ and $u_{xx}$ are color-coded black, blue, and red, respectively (same below).
As can be seen, the two types of lines are very close for all $N$ considered, indicating that $\alpha_{\rm R, M+}$ and $\beta_{\rm R, M}$ are good estimates for $\alpha_{\rm R, O}$ and $\beta_{\rm R, O}$, which also follows from comparing their values in Table~\ref{table_results_using_MS_plus_for_1d_pois_u_exp_m_x_m_0p5_square_dealii}, second row.

\begin{figure}[!ht]
   \subfloat[Type 2\label{mesh_type_2}]{\includegraphics[width=0.35\linewidth]{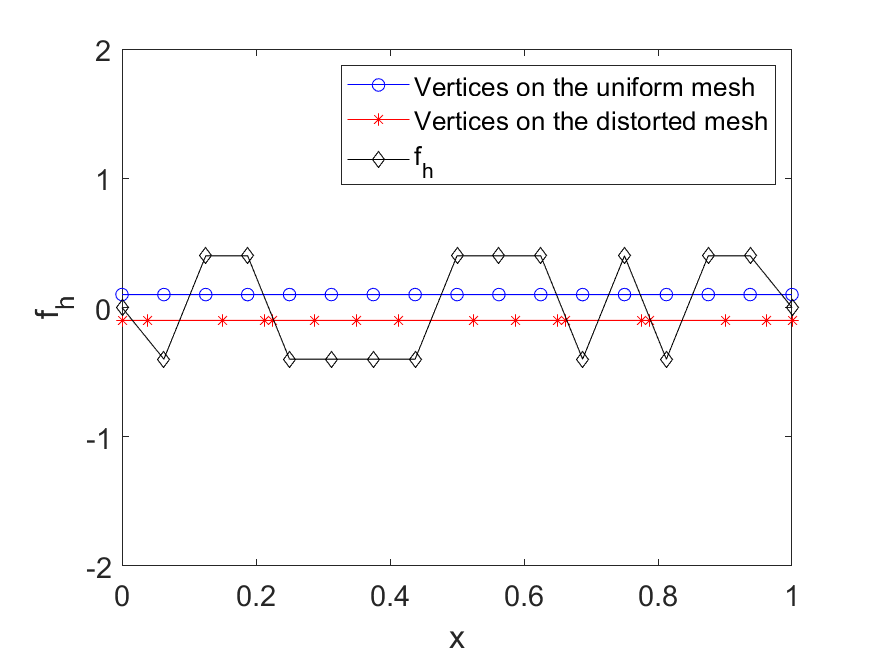}}
   \hspace{-1.5em}
   \subfloat[Type 3\label{mesh_type_3}]{\includegraphics[width=0.35\linewidth]{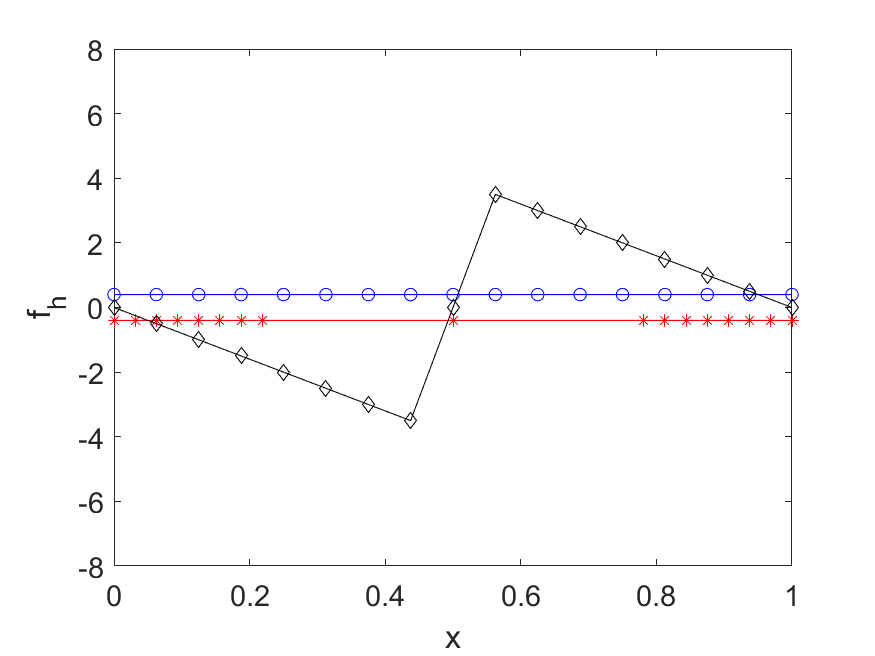}}
   \hspace{-1.5em}
   \subfloat[Type 4\label{mesh_type_4}]{\includegraphics[width=0.35\linewidth]{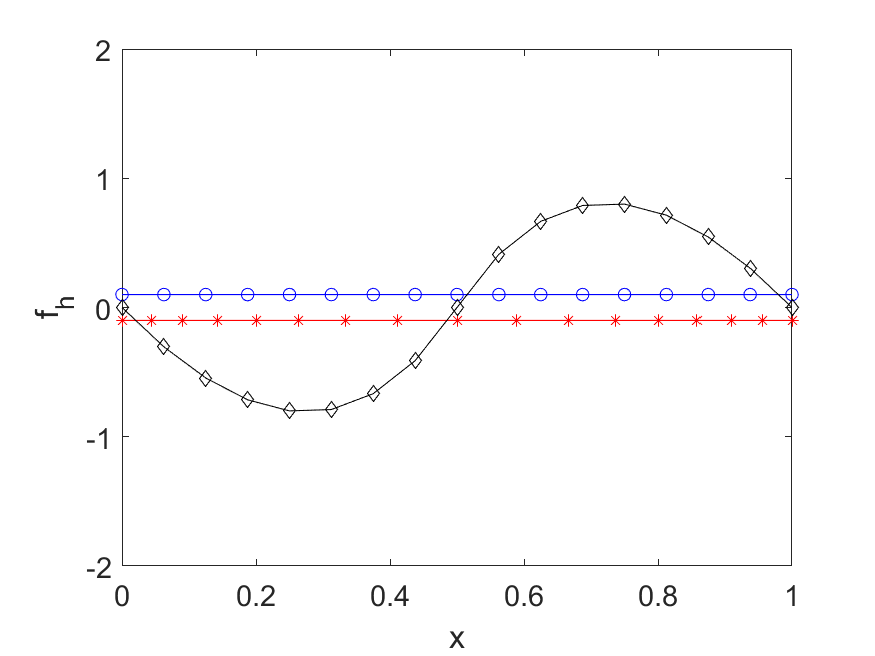}}     
   \caption{Vertex distribution for the distorted meshes.}
   \label{mesh_type_all}
\end{figure}

Next, we consider the influence of mesh distortion on $E_{\rm R}$.
Denoting the coordinate of a vertex on the equidistant mesh by $x_0$ and that on a distorted mesh by $x_1$, the distortion degree of a vertex is defined as
\begin{equation*}
  f_h = \frac{h_{\rm d}}{h_{\rm 0}},
\end{equation*}
where $h_{\rm d} = x_1 - x_0$ denotes the distorted distance of a vertex, and $h_{\rm 0}$ the grid size of the equidistant mesh, which is a constant. Obviously, $h_{\rm d}$ is positive if a vertex is moved to the right and negative if a vertex is moved to the left.
Three kinds of distorted meshes are considered, which are denoted by Mesh Type 2--4, respectively. For Mesh Type 2, $x_0$ is distorted randomly with $\mid f_h \mid = 0.4$. For Mesh Type 3 and Type 4, vertices are symmetric about $x = 0.5$, and when $x_0 < 0.5$, $x_0$ is distorted according to $x_1 = 0.5 x_0$ for Mesh Type 3 and according to $x_1 = \frac{x_0}{1.5-x_0}$ for Mesh Type 4. The resulting vertex distribution for the refinement level $R=4$ can be found in Fig.~\ref{mesh_type_2}--Fig.~\ref{mesh_type_4}, respectively. In these figures, the vertex distribution on the equidistant mesh and parameter $f_h$ are also shown for comparison.

For Mesh Type 2--4, the resulting coefficients that determine the round-off error $E_{\rm R}$ are shown in Table~\ref{table_results_using_MS_plus_for_1d_pois_u_exp_m_x_m_0p5_square_dealii}, lines 4--6.
The round-off error line represented by $\alpha_{\rm R, M+}$ and $\beta_{\rm R, M}$ is compared to that represented by $\alpha_{\rm R, O}$ and $\beta_{\rm R, O}$ in Figs.~\ref{manufacturing_alpha_R_beta_R_0_pois_Type_2}--\ref{manufacturing_alpha_R_beta_R_0_pois_Type_4}, respectively. 
It turns out that $\alpha_{\rm R, M+}$ and $\beta_{\rm R, M}$ also give a good estimate of $\alpha_{\rm R, O}$ and $\beta_{\rm R, O}$ for non-equidistant grids.

\paragraph{Diffusion problems}			\label{section_1d_benchmark_problems_diff}

For diffusion problems, we consider both $D = 1+x$ and $D = e^{-(x-0.5)^2}$. 
Only Mesh Type 1 is studied.

\begin{table}[!ht]
\caption{The coefficients $\alpha_{\rm R}$ and $\beta_{\rm R}$ obtained using the method of MS+ and OS for the 1D diffusion problem with $u = e^{-(x-0.5)^2}$.}
\centering
\scriptsize
\begin{tabular}{c|c|c|c|c|c|c}          \hline
 \multirow{2}{*}{$D$} & \multirow{2}{*}{Variable} & \multicolumn{3}{c|}{MS+} & \multicolumn{2}{c}{OS} \\ \cline{3-7}
   &  & $\alpha_{\rm R,M}$ & $\beta_{\rm R,M}$ & $\alpha_{\rm R, M+}$ & $\alpha_{\rm R, O}$ & $\beta_{\rm R, O}$ \\ \hline
 $1+x$ & \makecell{$u$\\ $u_x$\\ $u_{xx}$} & \makecell{1.0e-18\\ 5.0e-18\\ 1.0e-16} & \makecell{1.8\\ 1.8\\ 2.0} & \makecell{8.4e-18\\ 4.2e-17\\ 8.4e-16} & \makecell{2e-18\\ 2e-18\\ 5e-16} & \makecell{2.0\\ 2.0\\ 2.0} \\ \hline
 $e^{-(x-0.5)^2}$ & \makecell{$u$\\ $u_x$\\ $u_{xx}$} & \makecell{2.0e-18\\ 1.0e-17\\ 1.0e-16} & \makecell{1.5\\ 1.5\\ 2.0}& \makecell{1.7e-17\\ 8.4e-17\\ 8.4e-16} & \makecell{2e-17\\ 1e-16\\ 5e-16} & \makecell{1.5\\ 1.5\\ 2.0} \\ \hline
\end{tabular}
\label{table_results_using_MS_plus_for_1d_diff_u_exp_m_x_m_0p5_square_dealii}
\end{table}

\begin{figure}[!ht]
\centering
	\subfloat[$D=1+x$\label{manufacturing_alpha_R_beta_R_1_diff_D_1px_Type_1}]{\includegraphics[width=0.33\linewidth]{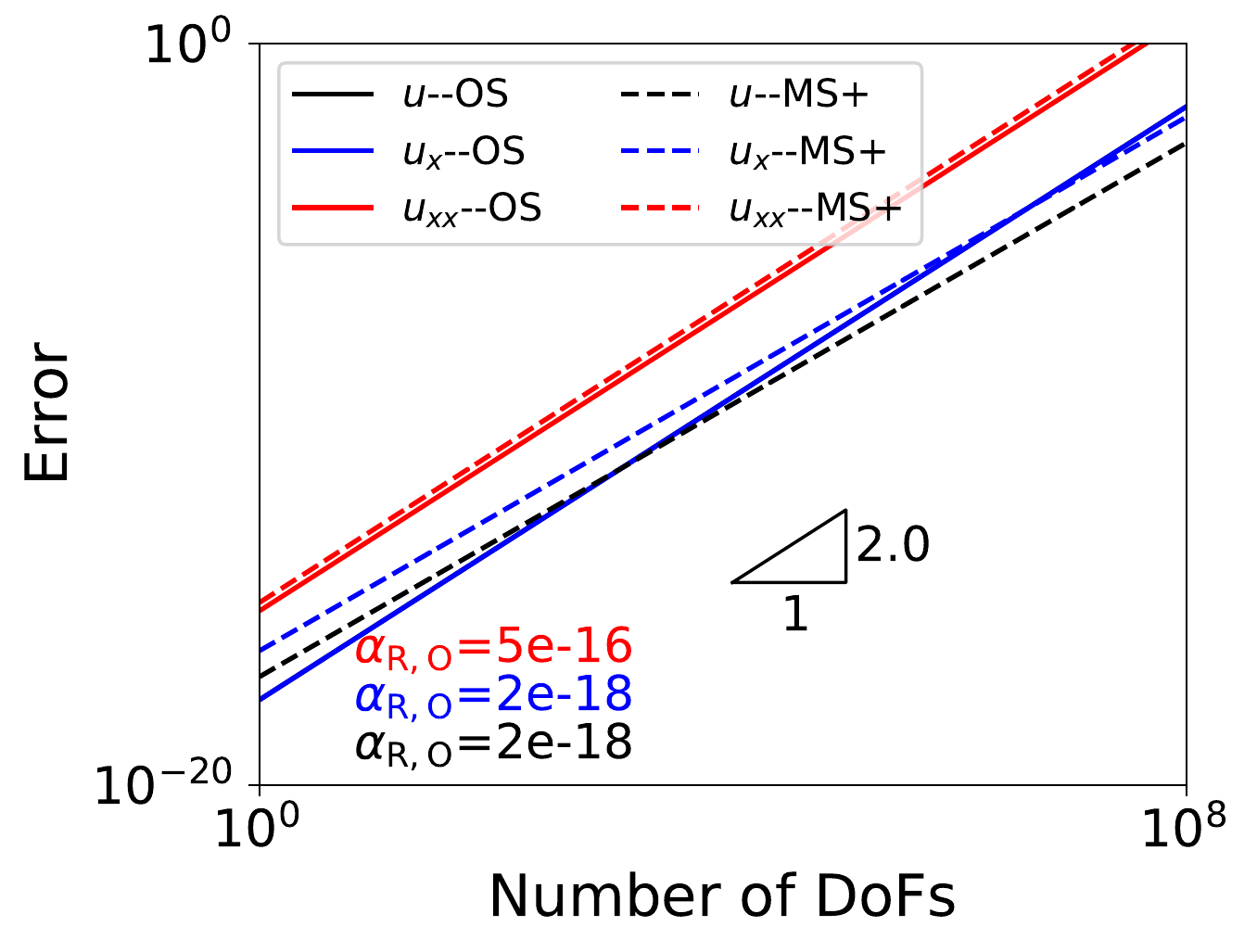}}
    \subfloat[$D=e^{-(x-0.5)^2}$\label{alpha_R_beta_R_1_diff_D_exp_m_x_m_0p5_square_Type_1}]{\includegraphics[width=0.33\linewidth]{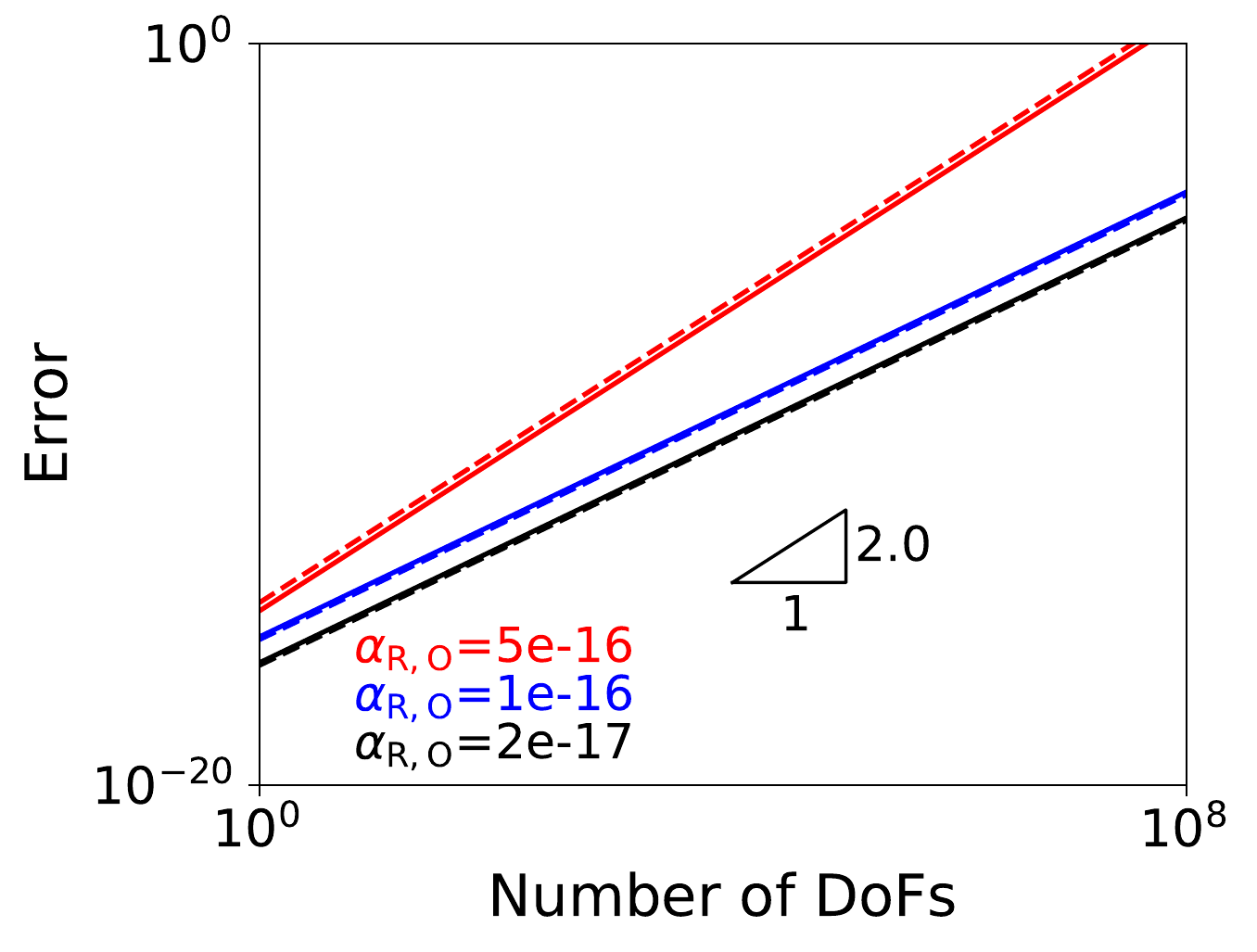}}
\caption{Comparison of the round-off error line represented by $\alpha_{\rm R, M+}$ and $\beta_{\rm R, M}$ to that represented by $\alpha_{\rm R, O}$ and $\beta_{\rm R, O}$ for the 1D diffusion problem with $u=e^{-(x-0.5)^2}$.}
\label{manufacturing_alpha_R_beta_R_1_diff}
\end{figure}

For both scenarios, the order of convergence of the truncation error is as expected; the resulting $\alpha_{\rm R, M}$, $\beta_{\rm R,M}$, $\alpha_{\rm R, M+}$, $\alpha_{\rm R, O}$, and $\beta_{\rm R, O}$ can be found in columns 3--6 of Table~\ref{table_results_using_MS_plus_for_1d_diff_u_exp_m_x_m_0p5_square_dealii}.
The round-off error lines obtained using $\alpha_{\rm R, M+}$ and $\beta_{\rm R, M}$ on the one hand and $\alpha_{\rm R, O}$ and $\beta_{\rm R, O}$ on the other hand are compared in Fig.~\ref{manufacturing_alpha_R_beta_R_1_diff}. As can be seen, the former gives a good approximation of the latter.

\begin{figure}[!ht]
	\subfloat[$r=1$\label{manufacturing_alpha_R_beta_R_2_helm_r_1}]{\includegraphics[width=0.33\linewidth]{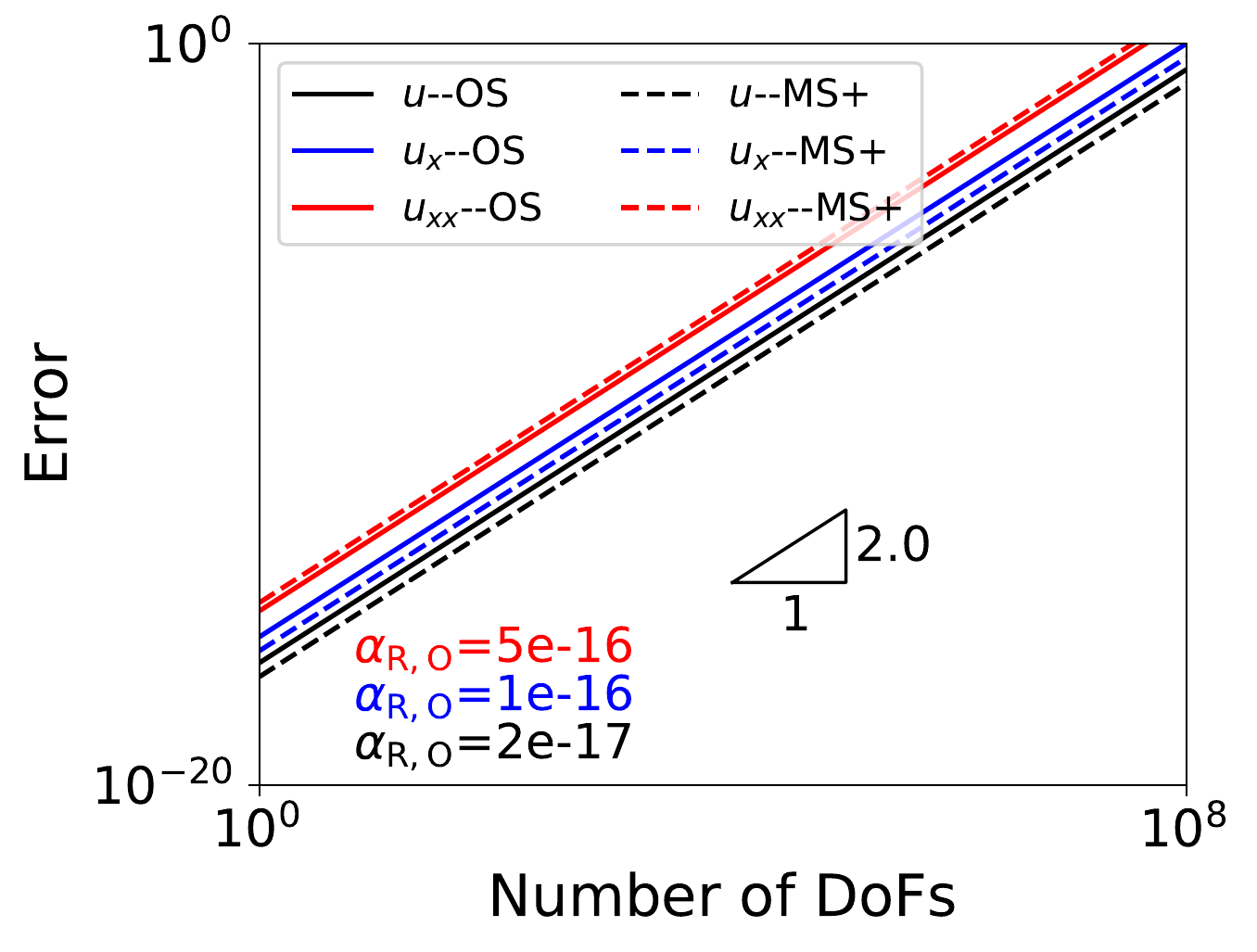}}
    \subfloat[$r=1+x$\label{alpha_R_beta_R_2_helm_r_1px}]{\includegraphics[width=0.33\linewidth]{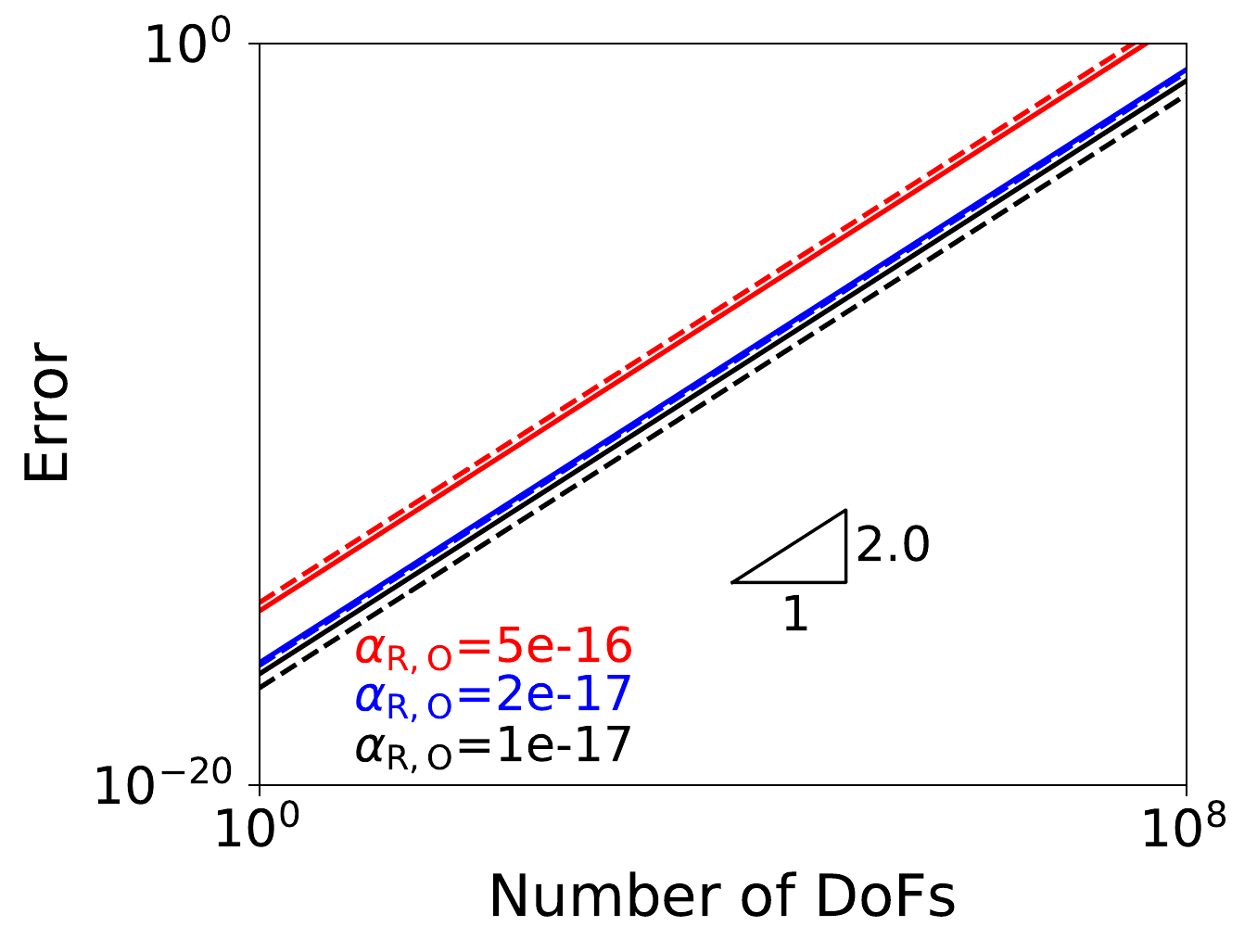}}
    \subfloat[$r=e^{-(x-0.5)^2}$\label{alpha_R_beta_R_2_helm_r_exp_m_x_m_0p5_square}]{\includegraphics[width=0.33\linewidth]{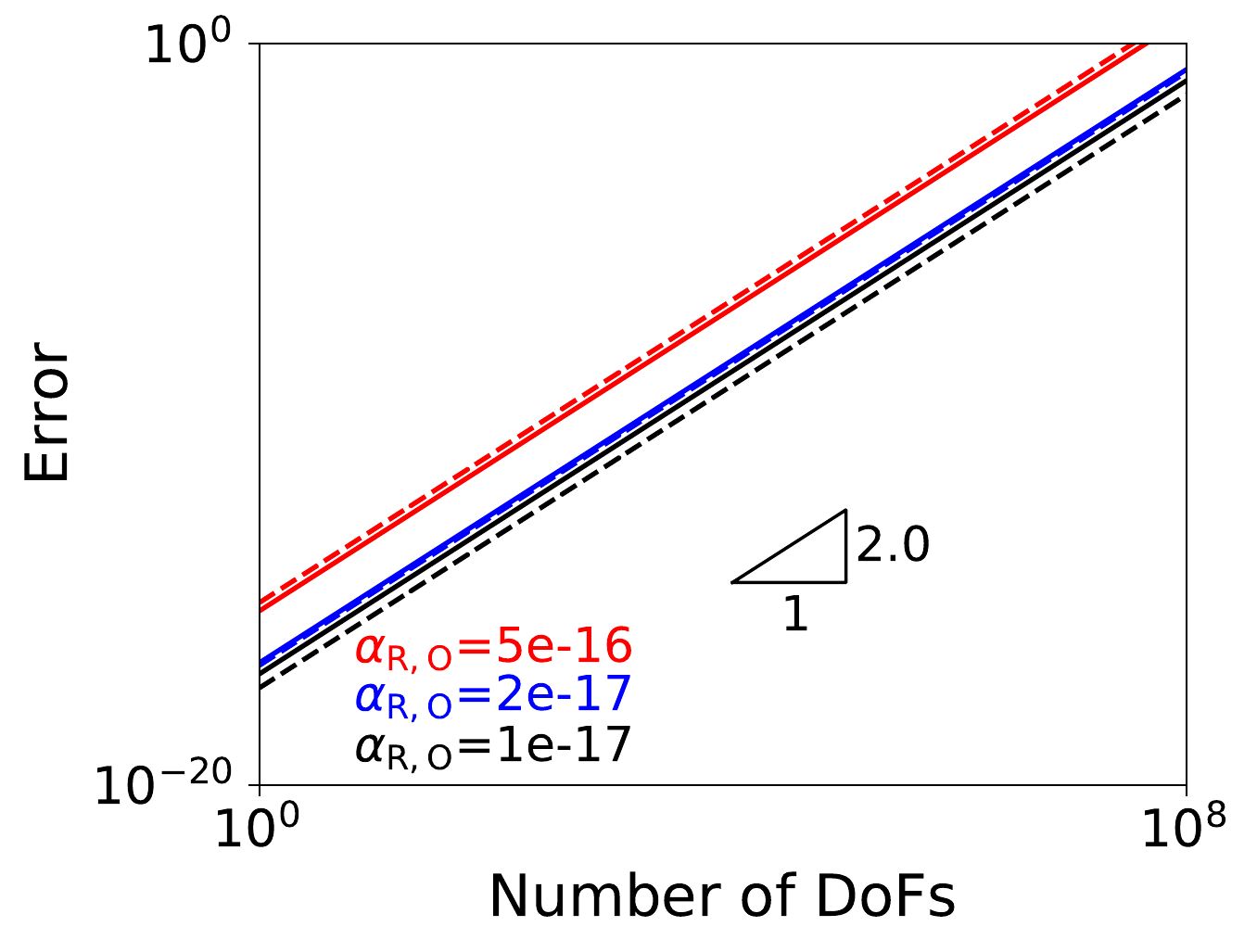}}
\caption{Comparison of the round-off error line represented by $\alpha_{\rm R, M+}$ and $\beta_{\rm R, M}$ to that represented by $\alpha_{\rm R, O}$ and $\beta_{\rm R, O}$ for the 1D Helmholtz problem with $u=e^{-(x-0.5)^2}$.}	
\label{manufacturing_alpha_R_beta_R_2_helm}
\end{figure}

\newpage
\paragraph{Helmholtz problems}				\label{section_1d_benchmark_problems_helm} 
For the Helmholtz problem, we consider $r= 1$, $r = 1+x$ and $r = e^{-(x-0.5)^2}$.
Only the uniform mesh is used.
For all the scenarios, the resulting $\alpha_{\rm R,M}$ and $\beta_{\rm R,M}$, $\alpha_{\rm R, M+}$, $\alpha_{\rm R, O}$, and $\beta_{\rm R, O}$ can be found in columns 3--7 of Table~\ref{table_results_using_MS_plus_for_1d_helm_u_exp_m_x_m_0p5_square_dealii}.
The round-off error lines represented by $\alpha_{\rm R, M+}$ and $\beta_{\rm R, M}$ are compared to that represented by $\alpha_{\rm R, O}$ and $\beta_{\rm R, O}$ in Fig.~\ref{manufacturing_alpha_R_beta_R_2_helm}. 
As can be seen, the two types of lines fit well. 
Therefore, the MS+ method is also suitable for predicting the dependency of the round-off error on the number of DoFs for general Helmholtz problems.

\begin{table}[!ht]
\caption{The coefficients $\alpha_{\rm R}$ and $\beta_{\rm R}$ obtained using the method of MS+ and OS for the 1D Helmholtz problem with $u = e^{-(x-0.5)^2}$.}
\centering
\scriptsize
\begin{tabular}{c|c|c|c|c|c|c}          \hline
 \multirow{2}{*}{$r$} & \multirow{2}{*}{Variable} & \multicolumn{3}{c|}{MS+} & \multicolumn{2}{c}{OS} \\ \cline{3-7}
   &  & $\alpha_{\rm R,M}$ & $\beta_{\rm R,M}$ & $\alpha_{\rm R, M+}$ & $\alpha_{\rm R, O}$ & $\beta_{\rm R, O}$ \\ \hline
 1 & \makecell{$u$\\ $u_x$\\ $u_{xx}$} & \makecell{1.0e-18\\ 5.0e-18\\ 1.0e-16} & \multirow{7}{*}{\makecell{2.0\\ 2.0\\ 2.0}} & \makecell{8.4e-18\\ 4.2e-17\\ 8.4e-16} & \makecell{2e-17\\ 1e-16\\ 5e-16} & \multirow{7}{*}{\makecell{2.0\\ 2.0\\ 2.0}} \\ \cline{1-3} \cline{5-6}
 $1+x$ & \makecell{$u$\\ $u_x$\\ $u_{xx}$} & \multirow{4}{*}{\makecell{5.0e-19\\ 2.0e-18\\ 1.0e-16}} &  & \multirow{4}{*}{\makecell{4.2e-18\\ 1.7e-17\\ 8.4e-16}} & \multirow{4}{*}{\makecell{1e-17\\ 2e-17\\ 5e-16}} &  \\  \cline{1-2}
 $e^{-(x-0.5)^2}$ & \makecell{$u$\\ $u_x$\\ $u_{xx}$} &  &  &  &  &  \\ \hline
\end{tabular}
\label{table_results_using_MS_plus_for_1d_helm_u_exp_m_x_m_0p5_square_dealii}
\end{table}

\newpage
\subsubsection{2D problems}                                \label{section_2d_benchmark_problems}
To assess the applicability of the MS+ method for 2D problems, we consider 2D Poisson and diffusion problems, with the solution given by $u = e^{-[(x-0.5)^2 + (y-0.5)^2]}$. 
The manufactured solution reads $u_{\rm M} = (x-0.5)^2+(x-0.5)(y-0.5)+(y-0.5)^2$.
$\|u_{\rm O}\|_2 \approx 0.86$ is about 4.1 times $\|u_{\rm M}\|_2 \approx 0.21$.

The diffusion matrix $D$ considered in the diffusion problem is given in the first column of Table~\ref{table_results_using_MS_plus_for_2d_problems_dealii}, third row. For all the scenarios, the resulting coefficients that determine the round-off error $E_{\rm R}$ are shown in Table~\ref{table_results_using_MS_plus_for_2d_problems_dealii}. As can be seen, $\alpha_{\rm R, M+}$ and $\beta_{\rm R, M}$ again approximate $\alpha_{\rm R, O}$ and $\beta_{\rm R, O}$ very well.

\begin{table}[!ht]
\centering
\caption{The coefficients $\alpha_{\rm R}$ and $\beta_{\rm R}$ obtained using the method of MS+ and OS for 2D problems with $u = e^{-[(x-0.5)^2 + (y-0.5)^2]}$.}
\scriptsize
\label{table_results_using_MS_plus_for_2d_problems_dealii}
\begin{tabular}{c|c|c|c|c|c|c}	\hline
 \multirow{2}{*}{$D$} & \multirow{2}{*}{Variable} & \multicolumn{3}{c|}{MS+} & \multicolumn{2}{c}{OS} \\ \cline{3-7}
   &  & $\alpha_{\rm R,M}$ & $\beta_{\rm R,M}$ & $\alpha_{\rm R, M+}$ & $\alpha_{\rm R, O}$ & $\beta_{\rm R, O}$ \\ \hline
    $\big(\begin{smallmatrix} 1 & 0 \\ 0 & 1 \end{smallmatrix} \big)$ & \makecell{$u$\\ $\nabla u$\\ $\mathbf{H} u$} & \makecell{1e-17\\ 2e-17\\ 5e-16} & \makecell{1.00\\ 1.00\\ 1.00} & \makecell{4.1e-17\\ 8.2e-17\\ 2.1e-15} & \makecell{5e-17\\ 2e-16\\ 2e-15} & \makecell{1.00\\ 1.00\\ 1.00} \\ \hline
    $\big(\begin{smallmatrix} 1+x+y & xy \\ xy & 1+x+y \end{smallmatrix} \big)$ & \makecell{$u$\\ $\nabla u$\\ $\mathbf{H} u$} & \makecell{1e-17\\ 2e-17\\ 2e-16} & \makecell{0.75\\ 0.75\\ 1.00} & \makecell{4.1e-17\\ 8.2e-17\\ 8.2e-16} & \makecell{2e-17\\ 1e-16\\ 2e-15} & \makecell{0.75\\ 0.75\\ 1.00} \\ \hline
\end{tabular}
\end{table}

\begin{figure}[!ht]
\centering
	\subfloat[Poisson\label{manufacturing_alpha_R_beta_R_0_pois_2d}]{\includegraphics[width=0.33\linewidth]{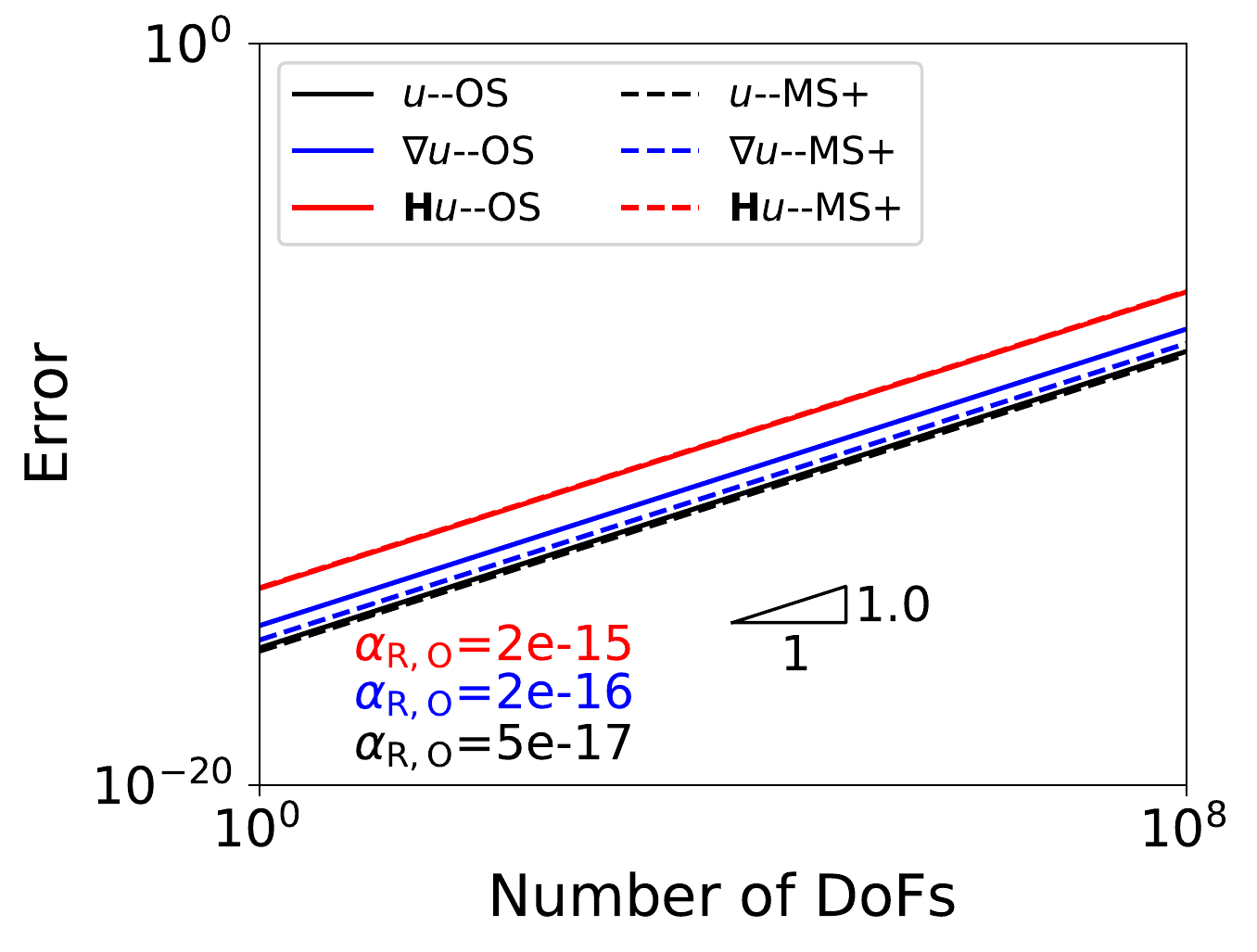}}
    \subfloat[Diffusion\label{alpha_R_beta_R_1_diff_2d}]{\includegraphics[width=0.33\linewidth]{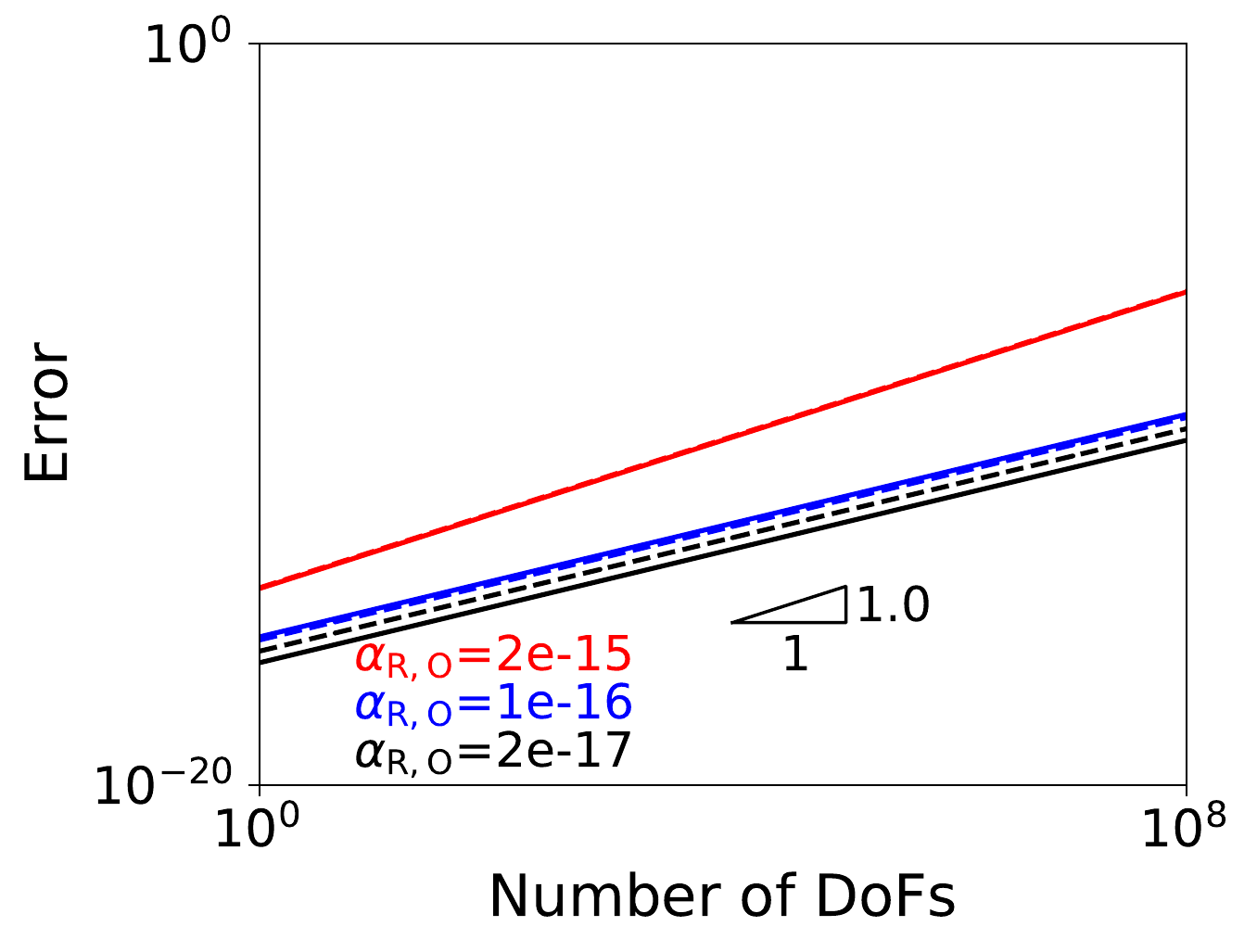}}
\caption{Comparison of the round-off error line represented by $\alpha_{\rm R, M+}$ and $\beta_{\rm R, M}$ to that represented by $\alpha_{\rm R, O}$ and $\beta_{\rm R, O}$ for the 2D problems with $u = e^{-[(x-0.5)^2 + (y-0.5)^2]}$.}
\label{manufacturing_alpha_R_beta_R_1_diff_2d}
\end{figure}


\section{Algorithm}					\label{section_algorithm}

In this section, we extend the a posteriori algorithm introduced in \cite{liu386balancing} to two-dimensional problems.
To keep this paper self-contained, we repeat the basic algorithmic steps.
The reader already familiar with our earlier work may want to skip this section, recalling that the essential modification of our algorithm is the use of a manufactured solution approach to derive $\alpha_{\rm R}$ and $\beta_{\rm R}$, cf. line 1 of algorithm 2.

In the algorithm, we define the following coefficients and use them in the steps given below.
\begin{itemize}
  \renewcommand\labelitemi{--}
  \item a minimal number of $h$-refinements before carrying out `\textit{PREDICTION}', denoted by $R_{\rm min}$, with the following default values for 1D problems:
  \begin{equation}
  \begin{aligned}
      R_{\rm min} &=
      \begin{cases*}
	9-p & for $p < 6$, \\
	4 & otherwise.
      \end{cases*}
  \end{aligned}
  \end{equation}
  Smaller values are chosen for 2D problems.
  We choose this parameter mainly because the error might increase, or decrease faster than the asymptotic order of convergence for coarse refinements, especially for lower-order elements.
  \item the allowed maximum $N$ : $10^8$, denoted by $N_{\rm max}$.  
  \item a stopping criterion $c_s = 0.1$ for seeking the $L_2$ norm of the dependent variable, denoted by $\|u\|_2$. When the difference of two adjacent $\|u_h\|_2$, evaluated by $\mid \frac{\|u _{h}\|_2 - \|u _{2h}\|_2}{\|u _{2h}\|} \mid$, is smaller than $c_s$, the iteration is stopped. We choose this parameter because the analytical solution does not exist for most practical problems.
  \item a relaxation coefficient $c_r$ for seeking the asymptotic order of convergence, with the following default values: 
    \begin{equation}
    \begin{aligned}
	c_r &=
	\begin{cases*}
	  0.9 & for $p < 4$, \\
	  0.7 & for 4 $\leqslant$ $p < 10$, \\
	  0.5 & otherwise.
	\end{cases*}
    \end{aligned}
    \end{equation}
\end{itemize}

The procedure of the algorithm consists of six steps, which are explained below:

\paragraph{Step-1} `\textit{INPUT}'. In this step, the items shown in the Table \ref{settings_algorithm} has to be provided by the user.

\begin{table}[!ht]
\small
\captionof{table}{Custom input of the algorithm.}
\label{settings_algorithm}
  \centering
  \begin{tabular}{l | L{10cm}}
    \toprule
    Type & Item  \\
    \midrule
    Problem & \tabitem the problem to be solved \\
     		& \tabitem variables of which the highest accuracy is of interest \\ \hline
    FEM     & \tabitem an ordered array of element degrees $\{p_{\min}, \ldots, p_{\rm max}\}$ \\
    \bottomrule
  \end{tabular}
\end{table}

\paragraph{Step-2} `\textit{NORMALIZATION}'. The function of this step is to find $\|u\|_2$, in which the element degree $p = 2$. The specific procedure can be found in Algorithm \ref{algo_scaling_factor}. 

\vspace{0.2cm}
\begin{algorithm}[H]
\caption{NORMALIZATION}
\label{algo_scaling_factor}
\While{$N<N_{\rm max}$}
{
    \eIf{$\left|\frac{\|u_{h}\|_2 - \|u_{2h}\|_2}{\|u_{2h}\|_2} \right| < c_s$}
    {
        $\|u\|_2$ $\gets$ $\|u_{h}\|_2$\;
        break\;
    }
    {
        $h$ $\gets$ $h/2$\;
        calculate $\|u_h\|_2$\;    
    }
}
\end{algorithm}

\paragraph{Step-3} `\textit{PARAMETERIZATION}'.
In this step, we determine $\alpha_{\rm R, M+}$ and $\beta_{\rm R, M}$. 
Let us remind the reader that this is the main modification of our algorithm introduced in \cite{liu386balancing} to make it work in 2D.
The procedure is summarized in Algorithm \ref{algo_parameterization} below.

\vspace{0.2cm}
\begin{algorithm}[H]
\caption{PARAMETERIZATION}
\label{algo_parameterization}
  Obtaining $\alpha_{\rm R, M}$, $\beta_{\rm R, M}$, $\| u_{\rm M} \|_2$, and $\| u_{\rm O} \|_2$\;
  Calculating $\alpha_{\rm R, M+}$\;
\end{algorithm}
                                
\paragraph{Step-4} `\textit{PREDICTION}'. This step finds $E_{\rm min}$ for each variable and $p$ of interest.
The procedure for carrying out this step can be found in Algorithm \ref{block_PREDICTION}.

\vspace{0.2cm}
\begin{algorithm}[H]
\caption{PREDICTION}
\label{block_PREDICTION}
    \While{$N<N_{\rm max}$ \textbf{\textup{and}} ${E_{h}}>E_{\rm R}$}
    {
        $q_h$ $\gets$ $\log _2 \left( {{E_{2h}}}/{{E_{h}}} \right)$\;
        \eIf
        {
            $q_h \geqslant \beta_{\rm T} \times c_r$
        }
        {
            $N_{\rm c} \gets N$\;
            $E_{\rm c} \gets {E_{h}}$\;
            $\alpha_{\rm T}$ $\gets$ ${E_{\rm c}}/{N_{\rm c}}^{- \beta_{\rm T}}$\;
            $N_{\rm opt} \gets \left( \frac{\alpha_{\rm T} \cdot \beta_{\rm T}}{\alpha _{\rm R} \cdot \beta_{\rm R}} \right)^{\frac{1}{\beta_{\rm R} + \beta_{\rm T}}}$\;
            $E_{\rm min} \gets \alpha_{\rm T} \cdot {N_{\rm opt}}^{- {\beta _{\rm T}}} + \alpha_{\rm R} \cdot {N_{\rm opt}}^{{\beta _{\rm R}}}$\;

        }
        {
            $h$ $\gets$ $h/2$\;
            calculate ${E_{h}}$\;
        }
	}    
\end{algorithm}

\paragraph{Step-5} `\textit{POSTPROCESS}'. In this step, $E_{\rm min} ^{\rm PRED+}$ is obtained by computing the result using $N_{\rm opt}$.

\paragraph{Step-6} `\textit{OUTPUT}'. In this step, we output $N_{\rm opt} ^{\rm PRED}$, $E_{\rm min} ^{\rm PRED+}$, and $T ^{\rm PRED+}$ obtained in \textit{Step-4} and \textit{Step-5}.

\section{Validation}                                \label{section_validation_PRED_plus}

In this section, we evaluate the accuracy and efficiency of the PRED+ method.
Denoting $N_{\rm opt}$, $E_{\rm min}$, and $T$ of the BF method by $N_{\rm opt} ^{\rm BF}$, $E_{\rm min} ^{\rm BF}$, and $T ^{\rm BF}$, respectively, and that of the PRED+ method by $N_{\rm opt} ^{\rm PRED+}$, $E_{\rm min} ^{\rm PRED+}$, and $T ^{\rm PRED+}$, respectively, the accuracy is validated by comparing $N_{\rm opt} ^{\rm PRED+}$ and $E_{\rm min} ^{\rm PRED+}$ with $N_{\rm opt} ^{\rm BF}$ and $E_{\rm min} ^{\rm BF}$, respectively, and the efficiency by comparing $T ^{\rm PRED+}$ with $T ^{\rm BF}$.
For the latter, to better understand the difference between $T ^{\rm PRED+}$ and $T ^{\rm BF}$, the number denoting the reduction of the CPU time in percentage
\begin{equation}
  \text{pct} = \frac{T ^{\rm BF} - T ^{\rm PRED+}}{T ^{\rm BF}} \cdot 100
\end{equation}
is introduced.

We investigate the 2D Poisson problem with $u=e^{-[(x-0.5)^2 + (y-0.5)^2]}$ introduced in Section \ref{section_2d_benchmark_problems}.
Both deal.\rom{2} and FEniCS are used, and the element degree ranges from 1 to 5.

\begin{figure}[!ht]
\centering
  \subfloat[deal.\rom{2}\label{py_2d_validation_N_opt_dealii}]{\includegraphics[width=0.33\linewidth]{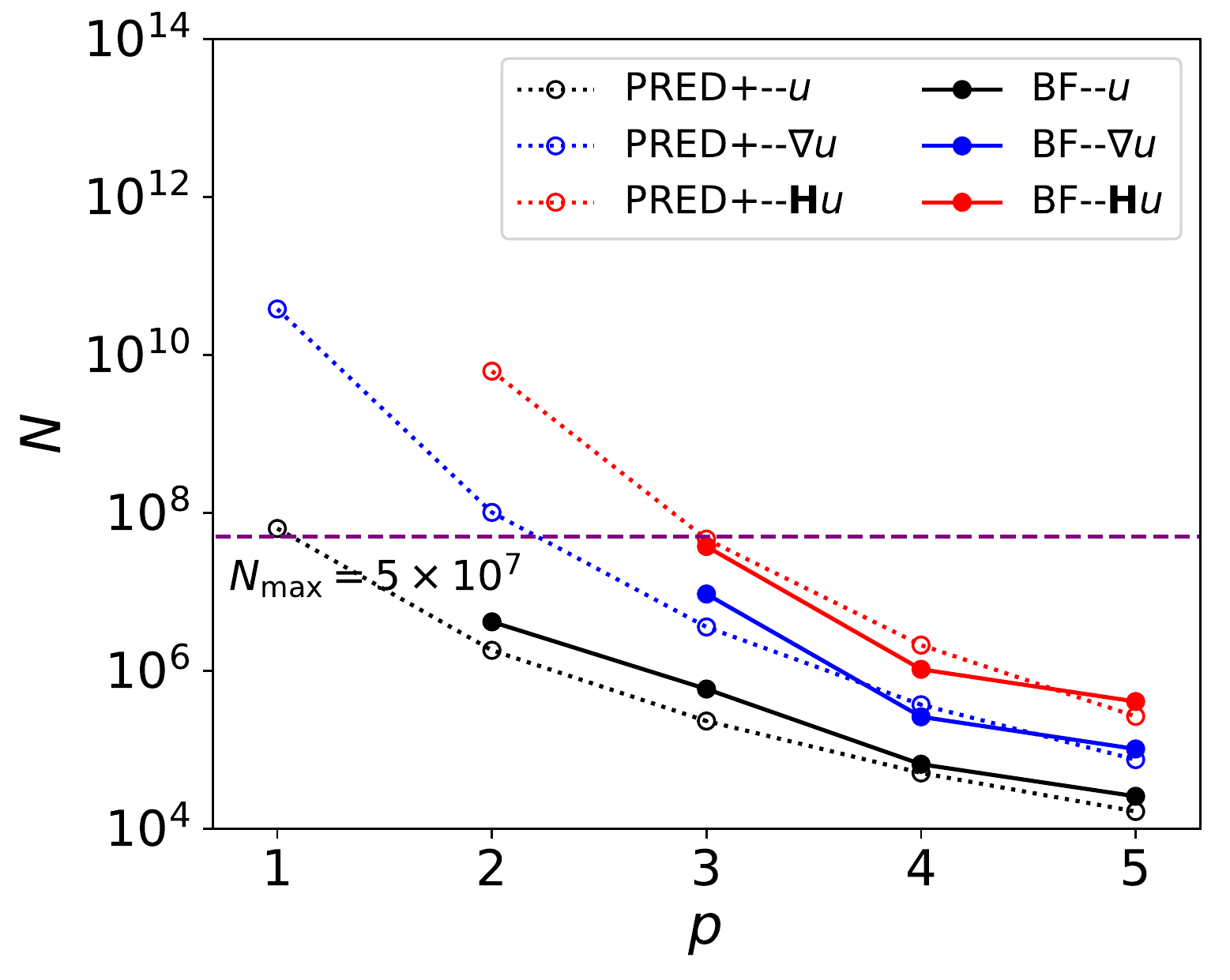}}   \quad
  \subfloat[FEniCS\label{py_2d_validation_N_opt_fenics}]{\includegraphics[width=0.33\linewidth]{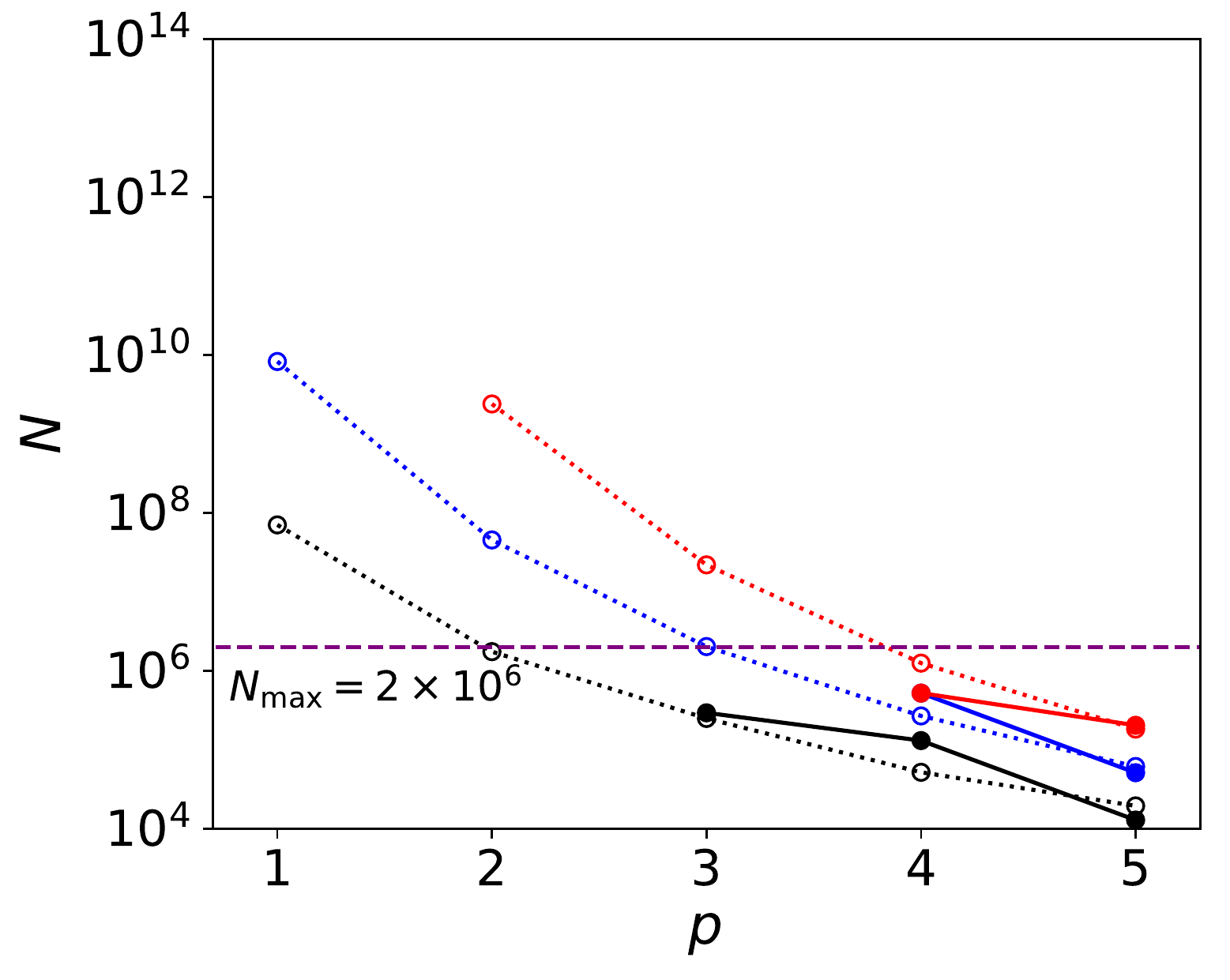}}
\caption{Comparison of $N_{\rm opt} ^{\rm PRED+}$ and $N_{\rm opt} ^{\rm BF}$ for the 2D Poisson problem with $u=e^{-[(x-0.5)^2 + (y-0.5)^2]}$.}
\label{py_2d_validation_N_opt}
\end{figure}

\begin{figure}[!ht]
\centering
  \subfloat[deal.\rom{2}]{\includegraphics[width=0.33\linewidth]{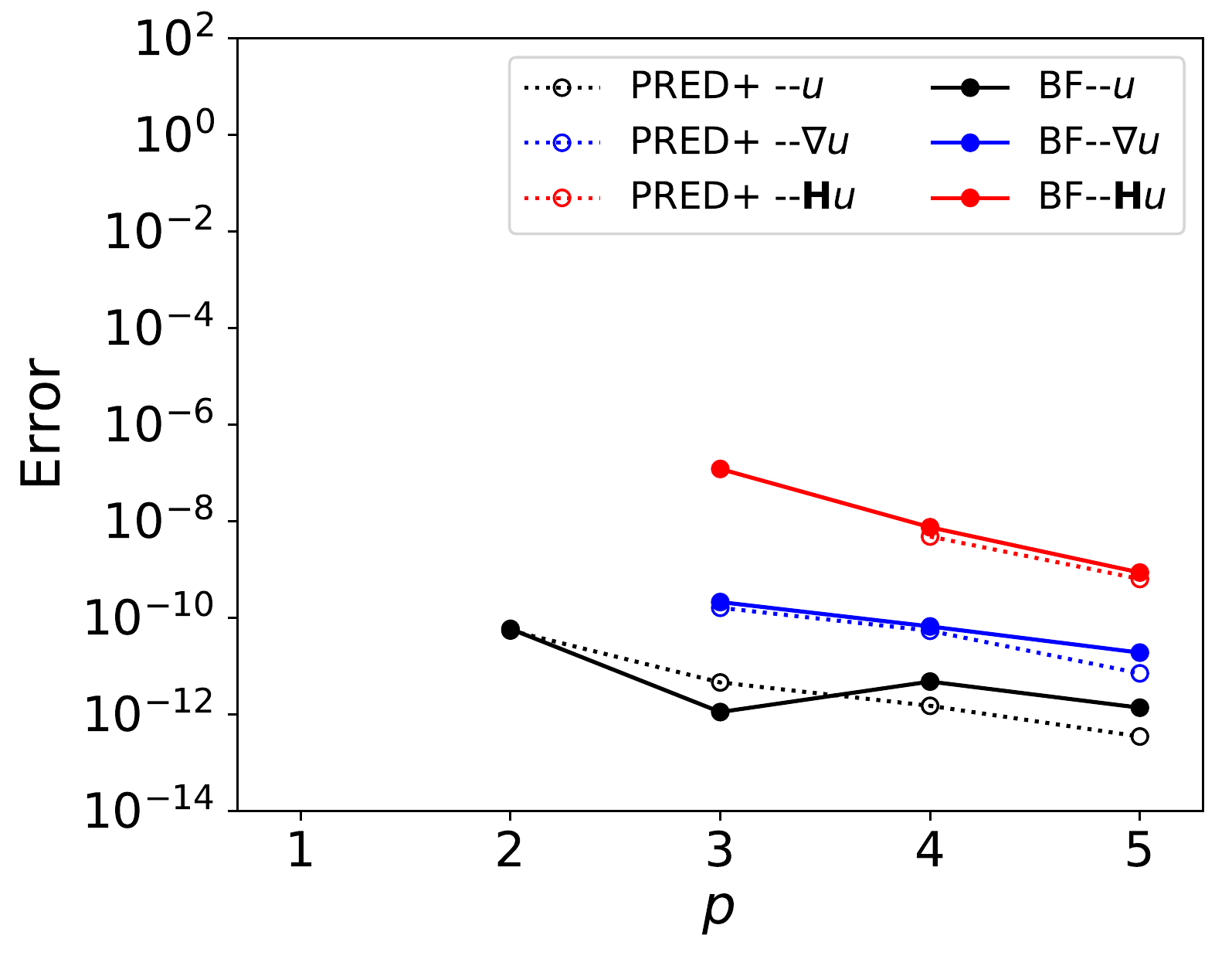}}   \quad
  \subfloat[FEniCS]{\includegraphics[width=0.33\linewidth]{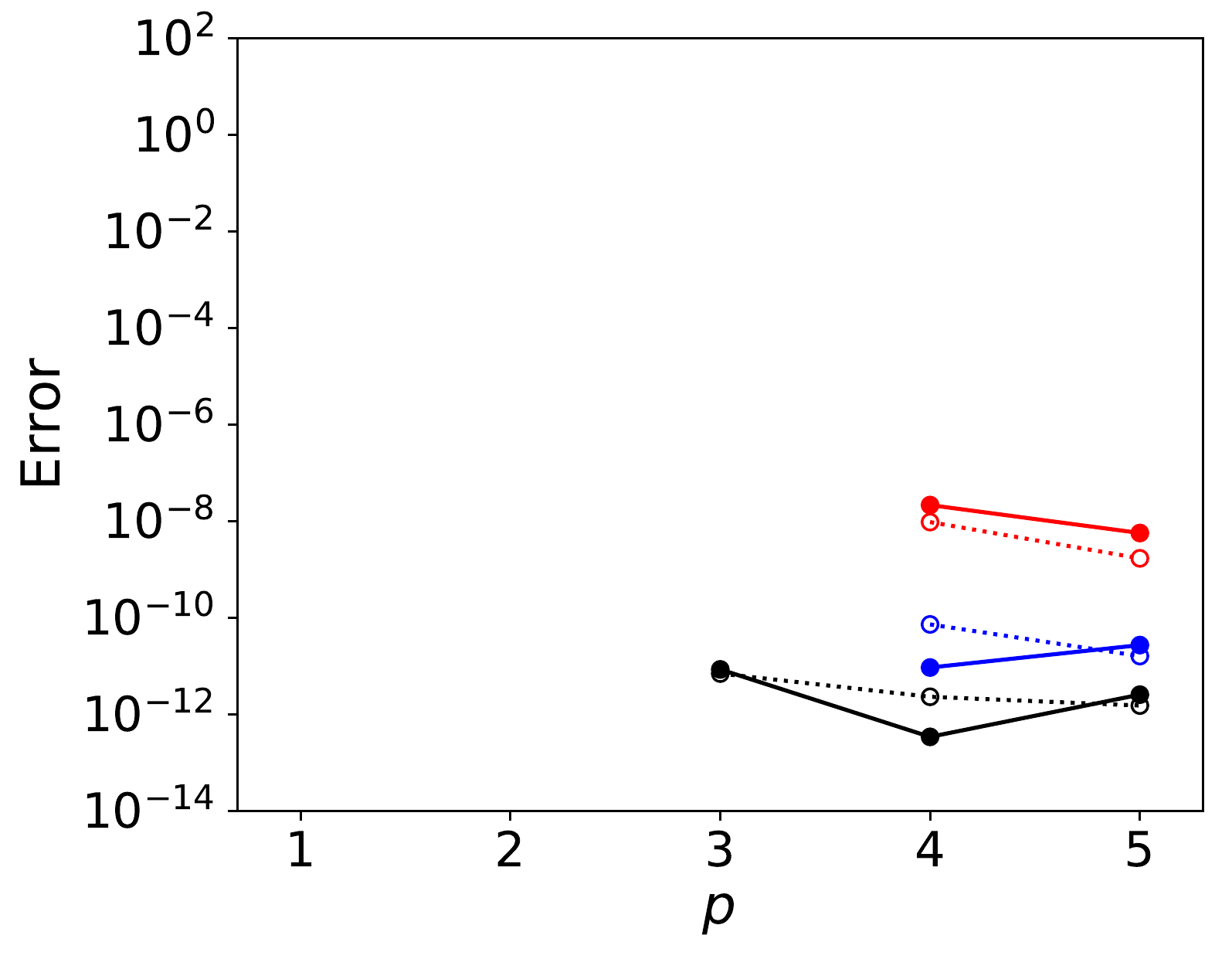}}  
\caption{Comparison of $E_{\rm min} ^{\rm PRED+}$ and $E_{\rm min} ^{\rm BF}$ for the 2D Poisson problem with $u=e^{-[(x-0.5)^2 + (y-0.5)^2]}$.}
\label{py_2d_validation_result_E_min}
\end{figure}

\begin{figure}[!ht]
\centering
  \subfloat[deal.\rom{2}]{\includegraphics[width=0.33\linewidth]{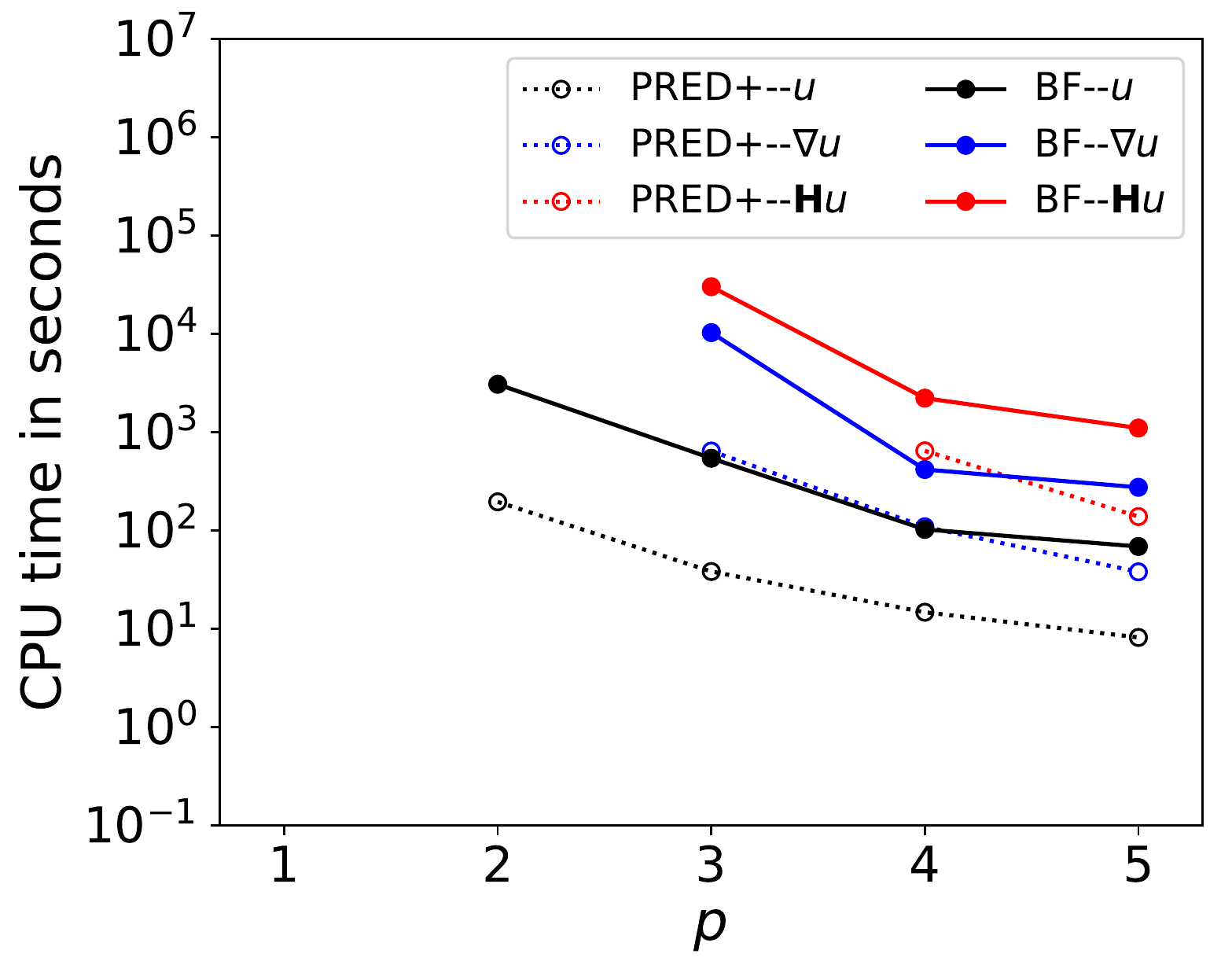}}   \quad
  \subfloat[FEniCS]{\includegraphics[width=0.33\linewidth]{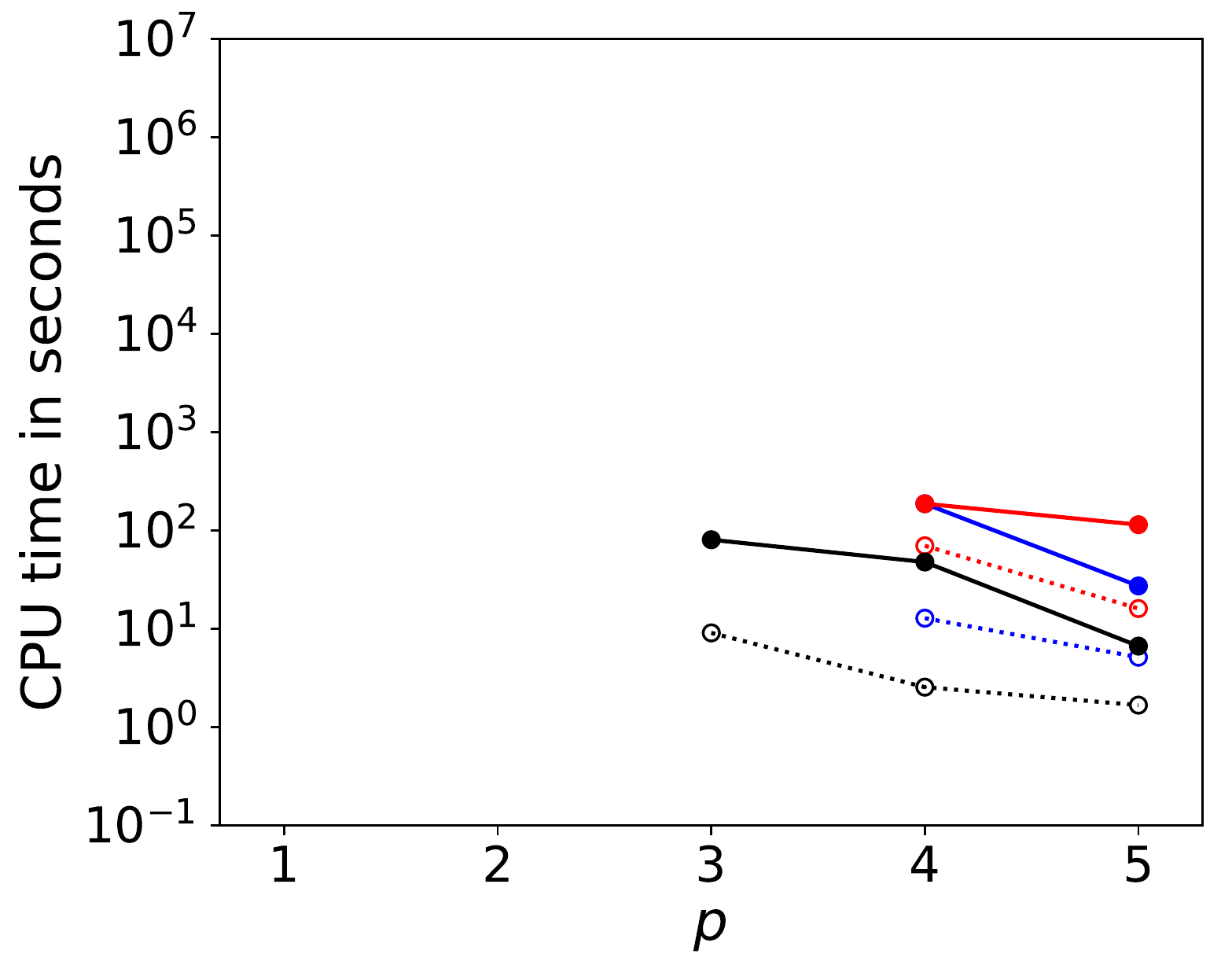}}  
\caption{Comparison of $T ^{\rm PRED+}$ and $T ^{\rm BF}$ for the 2D Poisson problem with $u=e^{-[(x-0.5)^2 + (y-0.5)^2]}$.}
\label{py_2d_validation_CPU_time_comparison}
\end{figure}

\begin{figure}[!ht]
\centering
  \subfloat[deal.\rom{2}]{\includegraphics[width=0.35\linewidth]{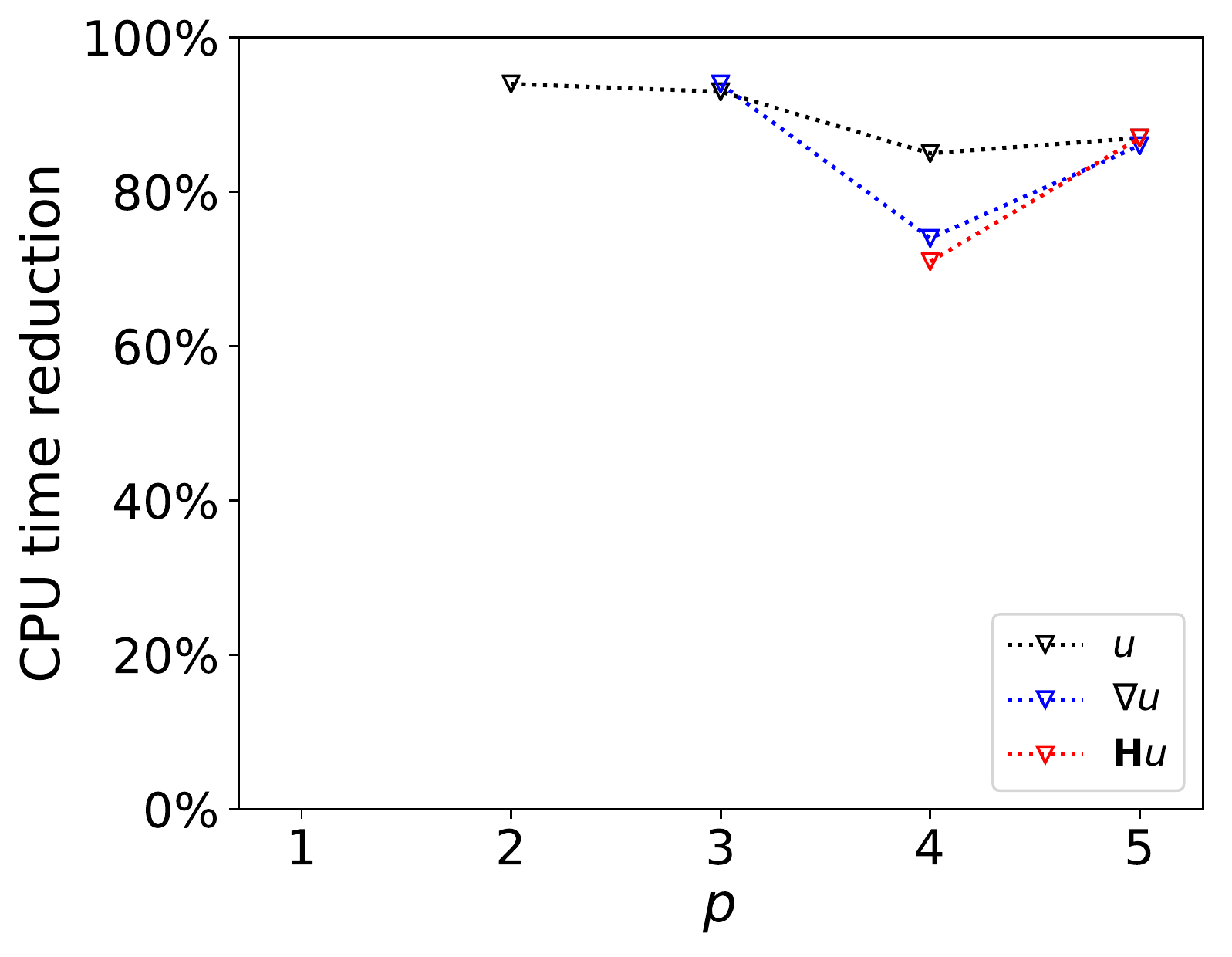}}   \quad
  \subfloat[FEniCS]{\includegraphics[width=0.35\linewidth]{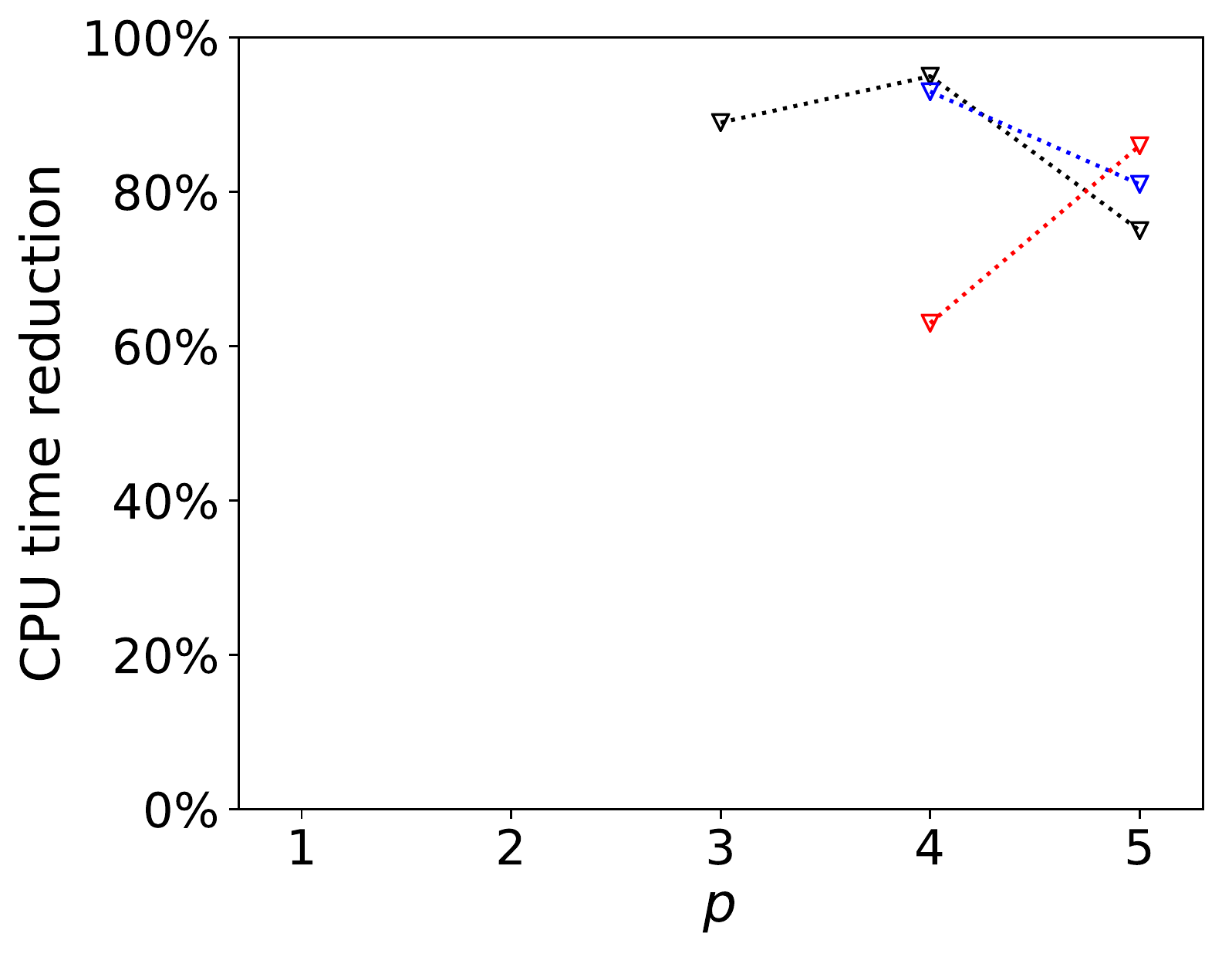}}  
\caption{CPU time reduction by the PRED+ method for obtaining different variables with the highest accuracy for the 2D Poisson problem with $u=e^{-[(x-0.5)^2 + (y-0.5)^2]}$.}
\label{py_2d_validation_CPU_time_reduction}
\end{figure}

Using the BF method, $N_{\rm opt} ^{\rm BF}$, $E_{\rm min} ^{\rm BF}$, and $T ^{\rm BF}$ are shown in Fig.~\ref{py_2d_validation_N_opt}--\ref{py_2d_validation_CPU_time_comparison}, respectively, and indicated by the solid dots connected by the solid line.
In Fig.~\ref{py_2d_validation_N_opt}, for some scenarios, the solution that corresponds to $N_{\rm opt} ^{\rm BF}$ DoFs cannot be obtained explicitly because the required number of DoFs is larger than $N_{\rm max}$, of which the value is shown by the purple dashed line. 
Using deal.\rom{2}, it concerns $p = 1$ for $u$, $p = 1$ and 2 for $\nabla u$, and $p = 2$ for $\mathbf{H} u$; using FEniCS, these are $p = 1,~2$ for $u$, $p = 1,~2$ and 3 for $\nabla u$, and $p = 2,~3$ for $\mathbf{H} u$.
Consequently, there is not corresponding data for $E_{\rm min} ^{\rm BF}$, $T ^{\rm BF}$, and pct in Fig.~\ref{py_2d_validation_result_E_min}--\ref{py_2d_validation_CPU_time_reduction}.

From the available data, it is found that $N_{\rm opt} ^{\rm BF}$, $E_{\rm min} ^{\rm BF}$, and $T ^{\rm BF}$ basically decrease for higher $p$ and increases for higher-order derivatives.
Therefore, the overall highest accuracy is controlled by the highest accuracy of the highest-order derivative when the variables required involve derivatives;
obtaining a high accuracy may be impossible using lower element degrees, while this is achievable using higher element degrees, for which the CPU time can be reduced.

Using the PRED+ method, $N_{\rm opt} ^{\rm PRED+}$, $E_{\rm min} ^{\rm PRED+}$, and $T ^{\rm PRED+}$ can be found in Fig.~\ref{py_2d_validation_N_opt}--\ref{py_2d_validation_CPU_time_comparison}, respectively. 
They are indicated by the open circles connected by the dotted line.
We are able to predict $N_{\rm opt}$ for all the scenarios.
For the scenarios with $N_{\rm opt} ^{\rm PRED+}$ larger than $N_{\rm max}$, the optimal solution cannot be obtained.
As a result, there is no corresponding data for $E_{\rm min} ^{\rm PRED+}$ and $T ^{\rm PRED+}$.

From the comparison for $N_{\rm opt}$, $E_{\rm min}$, and $T$ using the BF method and the PRED+ method, $N_{\rm opt} ^{\rm PRED+}$ and $E_{\rm min} ^{\rm PRED+}$ are very close to $N_{\rm opt} ^{\rm BF}$ and $E_{\rm min} ^{\rm BF}$, respectively.
However, the runtime $T ^{\rm PRED+}$ is much less than $T ^{\rm BF}$, see Fig.~\ref{py_2d_validation_CPU_time_comparison}. 
The CPU time reduction by the PRED+ method is about 60\% to 90\%, see Fig.~\ref{py_2d_validation_CPU_time_reduction}.
In summary, using the PRED+ method, we are able to determine if the accuracy required can be satisfied efficiently, and the CPU time for computing the result with the highest accuracy is reduced a lot.

\newpage

\section{Conclusion}					\label{paragraph_on_conclusion}

In this paper, we investigated the dependence of the round-off error on the number of DoFs when solving different 1D and 2D second-order boundary value problems using the standard FEM.
According to our findings, the round-off error increases according to a power-law function of the number of DoFs, and hence can be represented using two coefficients. 
A manufactured solution approach is proposed to determine the expression of the round-off error efficiently and accurately.
Using the expression obtained above, we extended our strategy in \cite{liu386balancing} for predicting and obtaining the highest achievable accuracy for 1D problems using deal.\rom{2} to 2D cases, and considered FEniCS as well. 
Using our strategy, the highest achievable accuracy obtained is very close to that obtained using sequential $h$-refinements, and the CPU time required can be reduced by 60\%$\sim$90\%.

\appendix

\section{Counting of the number of support points}                       \label{section_counting_number_support_points}

We first introduce the properties of the support points and then illustrate the counting of the number of support points.

\subsection{Properties of support points}                           \label{properties_support_points}

The support points are on various geometric properties, such as vertices, lines, and quads.
On one element, the number of geometric properties is shown in rows 2--3 of Table~\ref{table_for_geometric_properties}.
On the whole domain with the refinement level $R$, the number of geometric properties is shown in rows 4--5 of the same table.
\begin{table}[!ht]
\centering
\scriptsize
\caption{Number of geometric properties related to counting the number of support points.}
\begin{tabular}{c|c|c|c|c} \hline
 & Mesh element & Number of vertices & Number of faces & Number of quads \\ \hline
 \multirow{2}{*}{On a element} & Quadrilateral & 4 & 4 & 1 \\ \cline{2-5}
 & Triangle & 3 & 3 &	1	\\ \hline
 \multirow{2}{*}{\makecell{On the whole\\ domain}} & Quadrilateral & $(2^R + 1)^2$ & $2 \cdot (2^R + 1) \cdot 2^R$ & $(2^R)^2$  \\ \cline{2-5}
 & Triangle & $(2^R + 1)^2 + (2^R)^2$ & $2 \cdot (2^R + 1) \cdot 2^R + 4 \cdot (2^R)^2$ & $4 \cdot (2^R)^2$ \\ \hline 
\end{tabular}
\label{table_for_geometric_properties}
\end{table}

\begin{figure}[!ht]
  \centering
  \subfloat[Quadrilateral\label{sketch_distribution_degrees_of_freedom_dealii}]{\includegraphics[width=0.15\linewidth]{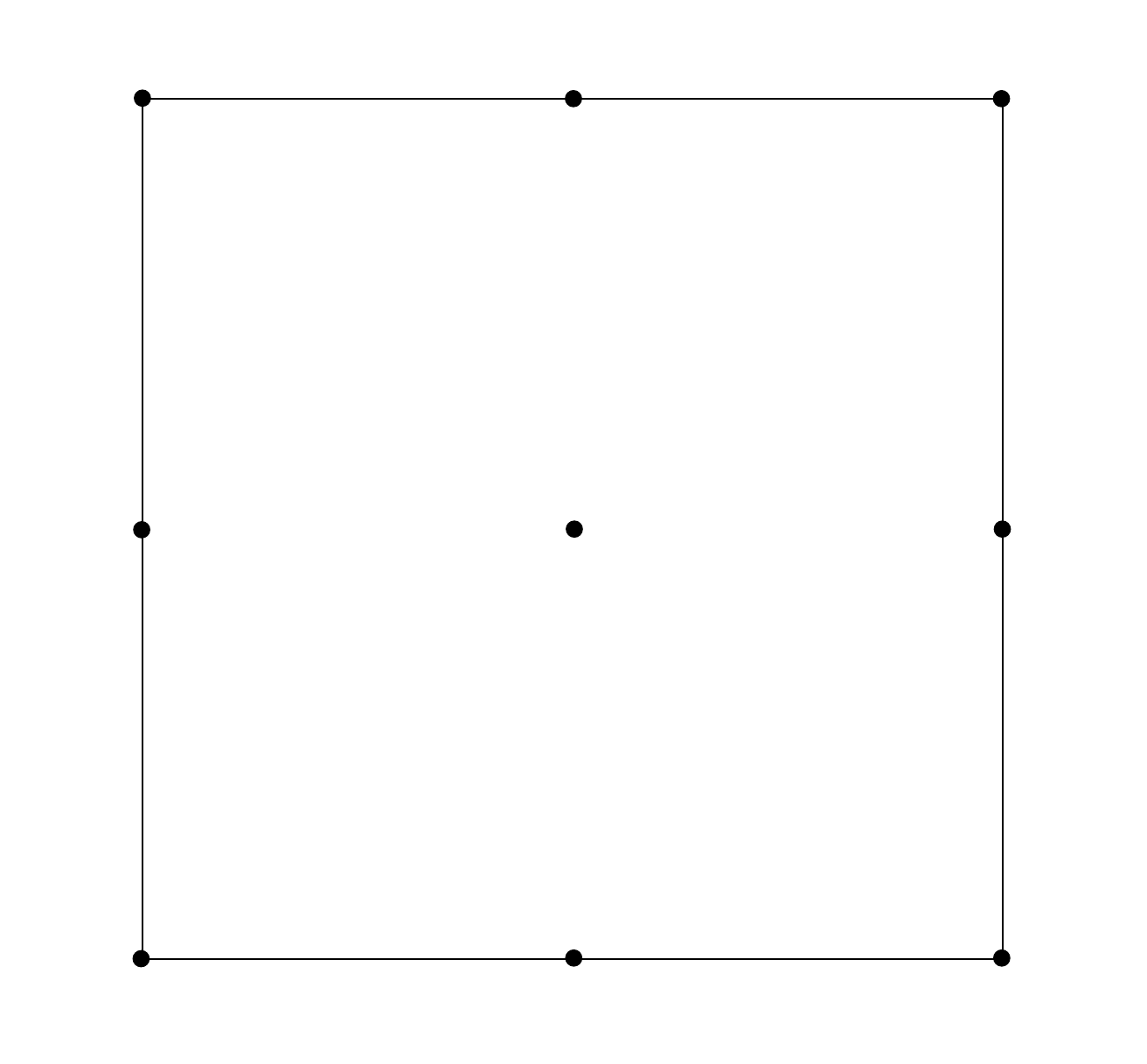}}
  \qquad
  \subfloat[Triangle\label{sketch_distribution_degrees_of_freedom_fenics}]{\includegraphics[width=0.15\linewidth]{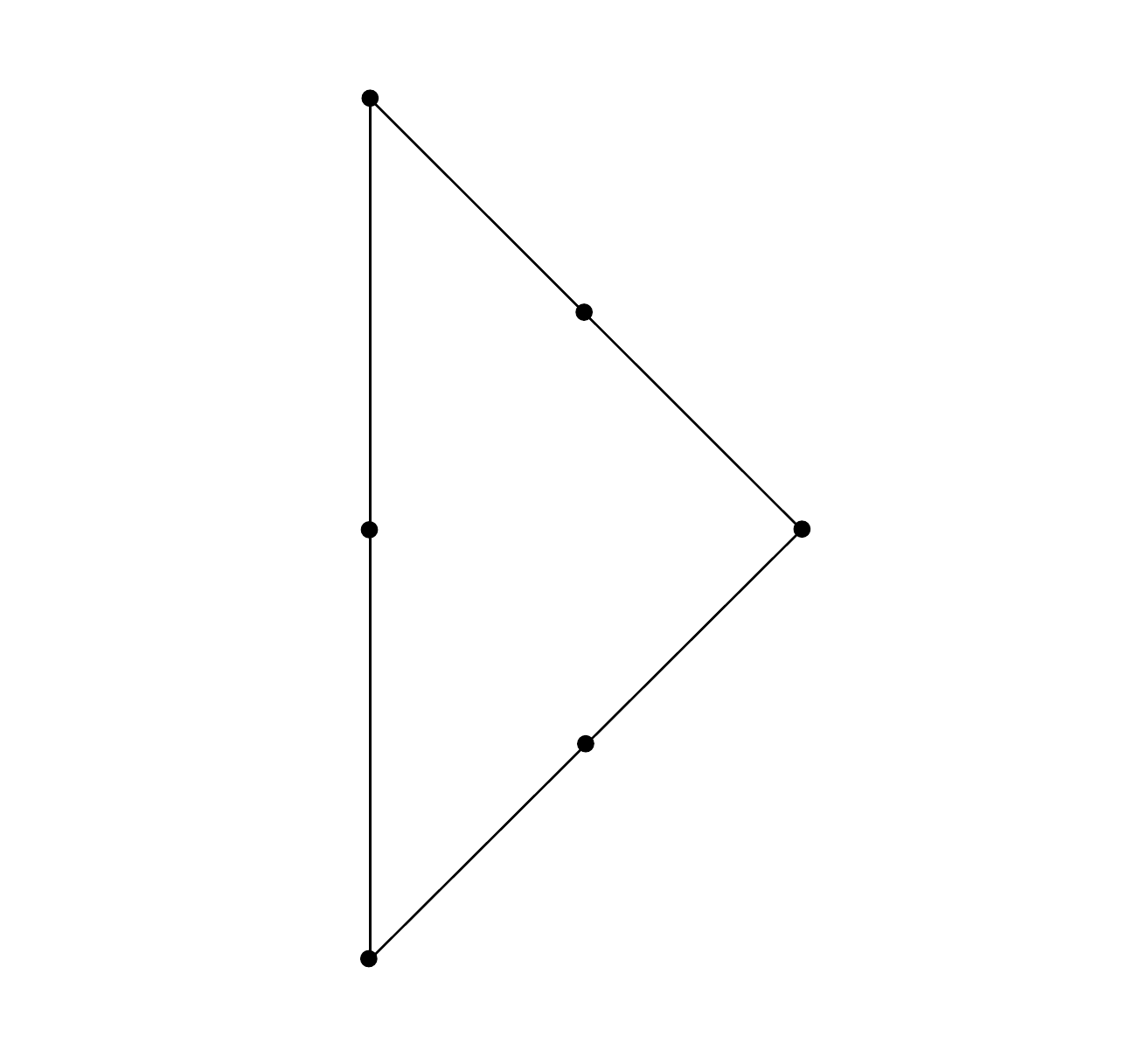}} 
   \caption{Location of support points (degrees of freedom) on an element when the element degree $p=2$.}
   \label{sketch_distribution_degrees_of_freedom}
\end{figure}

The number of support points on each geometric property is summarized in Table~\ref{table_for_number_of_support_points_on_geometric_properties}.
For example, when $p=2$, the location of support points can be found in Fig.~\ref{sketch_distribution_degrees_of_freedom}.
Using the numbers in Table~\ref{table_for_geometric_properties} and Table~\ref{table_for_number_of_support_points_on_geometric_properties}, we are able to count the number of support points over the whole domain.
\begin{table}[!ht]
\centering
\scriptsize
\caption{Number of support points on each geometric property of interest.}
\begin{tabular}{c|c|c|c} \hline
 FEM package & One vertex & One face & One quad \\ \hline
 deal.\rom{2} & 1 & $p-1$ & $(p-1)^2$ \\ \hline
 FEniCS & 1 & $p-1$ & $\frac{(p-1)(p-2)}{2}$	\\ \hline
\end{tabular}
\label{table_for_number_of_support_points_on_geometric_properties}
\end{table}

\subsection{Counting of the number of support points}

We show the number of support points on an element first and next the number of support points over the whole domain.
On one element, the number of support points reads
\begin{subequations}
\begin{align*}
  m_c &= 4 + 4 \cdot (p-1) + (p-1)^2	\\
      &= (p+1)^2				\numberthis
\end{align*}
using quadrilaterals and 
\begin{align}
  m_c = 3 + 3 \cdot (p-1) + \frac{(p-1)(p-2)}{2}
\end{align}
using triangles~\cite{arnold2014periodic}.
\end{subequations}
On the whole domain, the number of support points reads
\begin{subequations}
\begin{align*}
  m &= \underbrace{(2^R + 1)^2}_{\rm On~vertices} + \underbrace{2 \cdot (2^R + 1) \cdot 2^R \cdot (p-1)}_{\rm On~faces} + \underbrace{(2^R)^2 \cdot (p-1)^2}_{\rm On~quads} \\
  &= (2^R \cdot p + 1)^2	\numberthis 		\label{number_of_DoFs_on_the_whole_domain_dealii}
\end{align*}		
using quadrilaterals and 
\begin{align}
  m = \underbrace{\left[ (2^R + 1)^2 + (2^R)^2 \right]}_{\rm On~vertices} + \underbrace{\left[ 2^R \cdot (2^R + 1) \cdot 2 + 4 \cdot (2^R)^2 \right] \cdot (p-1)}_{\rm On~faces} + \underbrace{4 \cdot (2^R)^2  \cdot \frac{(p-2)(p-1)}{2} }_{\rm On~quads}				\label{number_of_DoFs_on_the_whole_domain_fenics}
\end{align}		
using triangles.
\label{number_of_DoFs_on_the_whole_domain}%
\end{subequations}


\section{Determination of \texorpdfstring{$\beta_{\rm T}$}{beta T}}				\label{section_proof_slope_ET}

This section is based on \ref{section_counting_number_support_points}.
Since the refinement level $R$ is relatively large when the analytical order of convergence is found, only keeping the terms containing $(2^R)^2$, we have 
\begin{subequations}
\begin{align}
  m \approx (2^R \cdot p)^2				\label{number_of_DoFs_Kh_approximation_dealii}
\end{align}
for Eq.(\ref{number_of_DoFs_on_the_whole_domain_dealii}), and 
\begin{align*}
  m &\approx \underbrace{\left[ 2 \cdot (2^R)^2 \right]}_{\rm On~vertices} + \underbrace{\left[ 6 \cdot (2^R)^2 \right] \cdot (p-1)}_{\rm On~faces} + \underbrace{2 \cdot (2^R)^2  \cdot (p^2 - 3 \cdot p + 2) }_{\rm On~quads} \\
  &= 2 \cdot (2^R)^2 \cdot p^2		\numberthis			\label{number_of_DoFs_Kh_approximation_fenics}
\end{align*}
for Eq.~(\ref{number_of_DoFs_on_the_whole_domain_fenics}).
\label{number_of_DoFs_Kh_approximation}%
\end{subequations}
\noindent From Eq.~(\ref{number_of_DoFs_Kh_approximation}), we have $\frac{N_{h/2}}{N_{h}} \approx 4$.

\bibliographystyle{unsrt}
\bibliography{bibfile_2d}

\end{document}